\documentclass{amsproc}

\usepackage[swedish,english]{babel}
\usepackage[utf8]{inputenc}
\usepackage[T1]{fontenc}
\usepackage{amsmath,amssymb,amsfonts}
\usepackage{url}
\usepackage[dvipsnames,usenames]{color}
\usepackage{array}
\usepackage[pdftex]{graphicx}
\usepackage{tikz}
\usepackage{accents}
\usepackage{tabu}

\DeclareFontFamily{U}{mathx}{\hyphenchar\font45}
\DeclareFontShape{U}{mathx}{m}{n}{
      <5> <6> <7> <8> <9> <10>
      <10.95> <12> <14.4> <17.28> <20.74> <24.88>
      mathx10
      }{}
\DeclareSymbolFont{mathx}{U}{mathx}{m}{n}
\DeclareFontSubstitution{U}{mathx}{m}{n}
\DeclareMathSymbol{\bigtimes}{1}{mathx}{"91}

\newcolumntype{L}{>{\displaystyle}l}
\newcolumntype{C}{>{\displaystyle}c}
\newcolumntype{R}{>{\displaystyle}r}

\renewcommand{\tfrac}{\genfrac{}{}{}1}

\newcommand{\R}{\mathbb R}

\newcommand{\Z}{\mathbb Z}
\newcommand{\C}{\mathbb C}

\newcommand{\D}{\mathbb D}
\renewcommand{\Im}{\mathrm{Im}}
\renewcommand{\Re}{\mathrm{Re}}
\def\E{{\mathrm{e}}}
 
\def\til{\widetilde}
\def\I{\mathfrak{i}}
\newcommand{\diff}{\mathrm{d}}
\renewcommand{\bar}{\overline}

\newcommand{\codim}{\mathrm{codim}}
\newcommand{\supp}{\mathrm{supp}}
\newcommand{\ntto}{\:\scriptsize{\xrightarrow{\wedge\:}}\:}

\newcommand{\ie}{\textit{i.e.}\/ } 
\newcommand{\eg}{\textit{e.g.}\/ } 
\newcommand{\cf}{\textit{cf.}\/ } 

\renewcommand{\vec}[1]{\accentset{\rightharpoonup}{#1}} 
\newcommand{\mat}[1]{#1} 

\allowdisplaybreaks

\theoremstyle{definition} 
\newtheorem{define}{Definition}[section]
\newtheorem{example}[define]{Example}
\newtheorem{remark}[define]{Remark}

\theoremstyle{plain} 

\newtheorem{thm}[define]{Theorem}
\newtheorem{prop}[define]{Proposition}
\newtheorem{coro}[define]{Corollary}

\numberwithin{equation}{section}

\begin{document}

\title[Geometric properties of measures]{Geometric properties of measures related to holomorphic functions having positive imaginary or real part}

\author{Annemarie Luger}
\address{Annemarie Luger, Department of Mathematics, Stockholm University, SE-106 91 Stockholm, Sweden}
\curraddr{}
\email{luger@math.su.se}
\thanks{The authors are supported by the Swedish Foundation for Strategic Research, grant nr. AM13-0011.}

\author{Mitja Nedic}
\address{Mitja Nedic, Department of Mathematics, Stockholm University, SE-106 91 Stockholm, Sweden, orc-id: 0000-0001-7867-5874}
\curraddr{}
\email{mitja@math.su.se}
\thanks{}


\date{2018-01-09} 

\begin{abstract}
In this paper, we study the properties of a certain class of Borel measures on $\R^n$ that arise in the integral representation of Herglotz-Nevanlinna functions. In particular, we find that restrictions to certain hyperplanes are of a surprisingly simple form and show that the supports of such measures can not lie within particular geometric regions, \eg strips with positive slope. Corresponding results are derived for measures on the unit poly-torus with vanishing mixed Fourier coefficients. These measures are closely related to  functions mapping the unit polydisk analytically into the right half-plane.
\end{abstract}

\maketitle

\section{Introduction}\label{sec:introduction}

When considering functions in one complex variable, those which map a complex half-plane or the unit disk into a half-plane play a  special role. They are very well studied and appear in many areas of applications. In particular, they are characterized via integral representations and there is a very intimate connection between the function and its representing measure, see \eg the classical article \cite{KacKrein1974}.

Generalizations of these classes to several complex variables have been considered, \eg in \cite{KoranyiPukanszky1963} for the unit polydisk, in \cite{Vladimirov1969} for the poly-upper half-plane and other tubular domains,  and have been investigated by several authors, \cf comments below.
It also has to be mentioned that these  functions in several complex variables  appear in applications as well, but are not yet utilized so much, due to the lack of understanding on the theoretical side, see \eg \cite{GoldenPapanicolaou1983,Milton2016}. There, the investigations are often restricted to two-component composite media, which requires functions in one variable only.

In the current text, we consider both functions that map the poly-upper half-plane $\C^{+n}$ or the unit polydisk $\D^{n}$ analytically into a half-plane. It is known that these functions can be characterized via integral representations of the form 
\begin{equation*}\label{eq:representation}
f(\vec{z})=L(\vec{z}) + \int_D K_n(\vec{z},\vec{t}) \diff\mu(\vec{t}). 
\end{equation*}
Here, $L$ denotes a linear term and $K_n$ a kernel function, both depending on the domain, while $D$ is the distinguished boundary of the domain and $\mu$ is a positive Borel measure on $D$, \cf Theorem \ref{thm:intRep_Nvar} and \cite[Theorem 1]{KoranyiPukanszky1963}. In the case of one variable, \ie $n=1$, all (reasonable) measures appear, whereas for several variables, not all measures are admissible. For functions mapping into the upper half-plane, the reason for that lies in the fact that the imaginary part of the kernel $K_n$ is then not non-negative, but may change sign. However, it can be decomposed as
\begin{equation}
\label{eq:kernel_intro}
\Im[K_n]=P_n+R_n,
\end{equation}
where $P_n$ denotes the respective Poisson kernel, which is positive, and $R_n$ denotes the remainder. It can than be shown that those measures which do appear as representing measures are exactly those which annihilate the remainder, \ie
\begin{equation}\label{eq:condition}
\int_D R_n(\vec{z},\vec{t}) \diff\mu(\vec{t}) = 0
\end{equation}
for all $z$, \cf Theorem \ref{thm:intRep_Nvar} and \cite[Theorem 1]{KoranyiPukanszky1963}. The drawback of these particular descriptions of the measures is that they are not so easy to check, which also makes it hard to construct examples. In the case of the unit polydisk, these representations and measures are discussed \eg in \cite{AglerHarlandRaphael2008,Forelli1981,KoranyiPukanszky1963,Rudin1969,Savchuk2006}. In particular, {examples of extremal measures are given in} \cite{McDonald1990}. For the poly-upper half-plane there are fewer results, \eg \cite{AglerHarlandRaphael2008,AizenbergDautov1976,Vladimirov1979} and, more recent, \cite{Savchuk2006}, as well as \cite{LugerNedic2016,LugerNedic2017,Nedic2017}.

We also want to mention that a subclass of these functions, namely the  Herglotz-Agler functions, are characterized via operator representations and so-called $\mu$-resolvents, \cf \cite{AglerEtal2012,AglerEtal2016,BallKaliuzhnyi-Verbovetskyi2015}.  However, these results do not imply the integral representation mentioned above and cannot be used for our current purpose. Moreover, there are many works on similar  functions defined on the ball, see \eg \cite{McCarthyPutinar2005} and, recently, \cite{Abu-GhanemAlpayColomboLewkowiczSabadini2018}. However, these results are not used here, since the ball and the polydisk (and, hence, also the poly-upper half-plane) are not biholomorphically equivalent.

In the present paper, we study the classes of representing measures, such that condition \eqref{eq:condition} is satisfied in the two cases considered. Particular focus lies on geometric properties of the support. However, contrary to many other texts concerning measures on the unit polydisk, we do not investigate the polydisk directly, but work, instead, first in the poly-upper half plane and then translate  the obtained results back to case of the polydisk.

In the poly-upper half-plane, holomorphic functions with non-negative imaginary part are called \emph{Herglotz-Nevanlinna functions}, \cf Definition \ref{def:HN_functions}, while their representing measures are called \emph{Nevanlinna measures}, \cf Definition \ref{def:Nevan_measure}.

First, we show that hyperplanes in $\mathbb R^n$ which are orthogonal to some coordinate axis play a special role for such measures. Namely, if the hyperplane is not a zero-set of the measure, then the restriction of the measure to that hyperplane has to be a constant multiple of the $(n-1)$-dimensional Lebesgue measure, \cf Theorem \ref{thm:affine_hyperplane_extraction}. In particular, the measure can be decomposed into one part supported on a coordinate orthogonal hyperplane and the remaining part, such that also the corresponding Herglotz-Nevanlinna function decomposes into two Herglotz-Nevanlinna functions, \cf Corollary \ref{coro:affine_hyperplane_extraction}. For other sets this procedure is not necessarily possible, since not every set can appear as the support of a Nevanlinna measure. We give several examples and discuss the situation for affine subspaces in detail. Moreover, we show that the support of a Nevanlinna measure cannot be confined within, roughly speaking, strips with positive slope, \cf Theorem \ref{thm:strips}. Applications of rational transformations yield even more results. In particular, we show that the support of a Nevanlinna measure cannot, for every coordinate, leave out a coordinate-orthogonal strip, \cf Theorem \ref{thm:cross}.

The structure of the paper is as follows. In the first part of this text, in Sections 2 and 3, we completely focus on the situation in the poly upper-half plane. In Section \ref{sec:prereq}, we review the integral representation theorem that lays the groundwork for our investigations. In Section \ref{sec:measures_plane}, we formally introduce the class of Nevanlinna measures on $\C^{+n}$ and present the main results of this paper, namely a detailed description of the form of these measures along coordinate-parallel affine subspaces of $\R^n$ as well as an investigation of the geometric properties of the support of such measures. The poly-torus will be discussed in Section \ref{sec:measures_polydisk}, where we investigate how the properties established for Nevanlinna measures in Section \ref{sec:measures_plane} relate back to measures on the unit poly-torus with vanishing mixed Fourier coefficients.

\section{Prerequisites}\label{sec:prereq}

We begin by recalling the following class of functions related to the poly-upper half-plane {$\C^{+n}:=\big\{z\in\C^n \,\big |\,\forall j=1,2,\ldots, n:  \Im[z_j]>0 \big\}$}.

\begin{define}\label{def:HN_functions}
A function $q\colon \C^{+n} \to \C$ is called a \emph{Herglotz-Nevanlinna function} if it is holomorphic and has non-negative imaginary part.
\end{define}

Our main tool in the study of Herglotz-Nevanlinna functions is the following characterization theorem \cite[Theorem 4.1 and {Theorem 5.1}]{LugerNedic2017}.

\begin{thm}\label{thm:intRep_Nvar}
A function $q\colon \C^{+n} \to \C$ is a Herglotz-Nevanlinna function if and only if $q$ can be written as
\begin{equation}\label{eq:intRep_Nvar}
q(\vec{z}) = a + \sum_{\ell=1}^nb_\ell z_\ell + \frac{1}{\pi^n}\int_{\R^n}K_n(\vec{z},\vec{t})\diff\mu(\vec{t}),
\end{equation}
where $a \in \R$, $\vec{b} \in [0,\infty)^n$, the kernel $K_n$ is defined for $\vec{z} \in \C^{+n}$ and $\vec{t} \in \R^n$ as
\begin{equation}\label{eq:kernel_Nvar}
K_n(\vec{z},\vec{t}) := \I\left(\frac{2}{(2\I)^n}\prod_{\ell=1}^n\left(\frac{1}{t_\ell-z_\ell}-\frac{1}{t_\ell+\I}\right)-\frac{1}{(2\I)^n}\prod_{\ell=1}^n\left(\frac{1}{t_\ell-\I}-\frac{1}{t_\ell+\I}\right)\right)
\end{equation}
and $\mu$ is a positive Borel measure on $\R^n$ satisfying the growth condition
\begin{equation}
\label{eq:measure_growth}
\int_{\R^n}\prod_{\ell=1}^n\frac{1}{1+t_\ell^2}\diff\mu(\vec{t}) < \infty
\end{equation}
and the Nevanlinna condition
\begin{equation}
\label{eq:measure_Nevan}
\int_{\R^n}\frac{1}{(t_{\ell_1} - z_{\ell_1})^2(t_{\ell_2} - \bar{z_{\ell_1}})^2}\prod_{\substack{j=1 \\ j \neq \ell_1,\ell_2}}^n\left(\frac{1}{t_j - z_j} - \frac{1}{t_j - \bar{z}_j}\right)\diff\mu(\vec{t}) = 0
\end{equation}
for all $\vec{z} \in \C^{+n}$ and all indices $\ell_1,\ell_2 \in \{1,2,\ldots,n\}$ with $\ell_1 < \ell_2$. Furthermore, for a given function $q$, the triple of representing parameters $(a,\vec{b},\mu)$ is unique.
\end{thm}

\begin{remark}
The Nevanlinna condition \eqref{eq:measure_Nevan} is, here, taken as one of the alternatives presented in \cite[Theorem 5.1]{LugerNedic2017} and is equivalent to the requirement that the measure annihilates the remainder term in \eqref{eq:kernel_intro}. 
\end{remark}

In the case $n = 1$, the above theorem reduces to the classical result attributed to Nevanlinna \cite{Cauer1932,Nevanlinna1922}. The case $n = 2$ was treated in \cite[Theorem 3.1]{LugerNedic2016} and an integral representation of the same form, but not as an "if and only if"-characterization with the accompanying conditions, appears also in \cite[{Section 17.4}]{Vladimirov1979}.

One consequence of Theorem \ref{thm:intRep_Nvar}, that will be of use several times later, is described by the corollary below, \cf \cite[Corollary 4.6]{LugerNedic2017}, where the symbol $\ntto$ denotes a non-tangential limit.

\begin{coro}\label{coro:propeties_of_HN_funcitons}
Let $n \geq 1$, let $q$ be a Herglotz-Nevanlinna function, let $p \in \R$ and let $j \in \{1,2,\ldots,n\}$. Then, there exists a non-negative number $c_j(p)$, such that
\begin{equation}\label{eq:c_constant_shifted}
\lim\limits_{z_j \ntto p}(p-z_j)\:q(\vec{z}) = c_j(p).
\end{equation}
In particular, the above limit is independent of the entries of the vector $\vec{z}$ at the non-$j$-th positions.
\end{coro}

\section{Properties of Nevanlinna measures}\label{sec:measures_plane}

Let us now formally introduce the class of measures we are going to study.

\begin{define}\label{def:Nevan_measure}
A positive Borel measure $\mu$ on $\R^n$ is called a \emph{Nevanlinna measure} if it satisfies both the growth condition \eqref{eq:measure_growth} and the Nevanlinna condition \eqref{eq:measure_Nevan}.
\end{define}

\begin{remark}
Note that a Nevanlinna measure is the representing measure of a whole family of Herglotz-Nevanlinna functions, as the linear part of representation \eqref{eq:intRep_Nvar} differs between different functions with the same representing measure.
\end{remark}

In the case $n = 1$, the class of Nevanlinna measures is merely the class of all positive Borel measure on $\R$ satisfying the growth condition $\int_\R(1+t^2)^{-1}\diff\mu(t) < \infty$. In higher dimensions, the appearance of the Nevanlinna condition \eqref{eq:measure_Nevan}, which conveniently reduces to an empty condition when $n = 1$, makes these measures much more interesting and involved.

In particular, in the case $n=2$, it was shown that the Nevanlinna condition \eqref{eq:measure_Nevan} implies  for non-trivial Borel measures that they cannot be finite as well as that points have zero mass, \ie $\mu(\{\vec{t}_0\}) = 0$ for any point $\vec{t}_0 \in \R^2$, see \cite[Propositions 4.3 and 4.4]{LugerNedic2016}. Among others, we will show that corresponding results hold for all $n\geq2$. The first will be formulated in Proposition \ref{prop:infinite_measure}, whereas the second turns out to be a special case of Theorem \ref{thm:affine_hyperplane_extraction}.

\subsection{Mass of affine subspaces}\label{subsec:subspaces}

We start with the following proposition on the measure of the whole space, which already marks a big difference between dimension $1$ and higher dimensions. 

\begin{prop}\label{prop:infinite_measure}
A non-trivial Nevanlinna measure cannot be finite.
\end{prop}

The proof goes along the same lines as for \cite[Proposition 4.3]{LugerNedic2016} and is, hence, omitted here. 

For a Nevanlinna measure $\mu$, we are now turning to its restrictions to hyperplanes in $\R^n$, which are orthogonal to some coordinate axis. More precisely, for a given index $j \in \{1,2,\ldots,n\}$ and a given point $p \in \R$, we denote the affine hyperplane 
$$H_j(p) := \{\vec{t} \in \R^n~|~t_j = p\} \subseteq \R^n.$$
A main result of this paper is the following description of the measure $\mu$ along such a hyperplane.

\begin{thm}\label{thm:affine_hyperplane_extraction}
Let $n \geq 2$ and let $\mu$ be a Nevanlinna measure. Take an index $j \in \{1,2,\ldots,n\}$ and any point $p \in \R$. Let $\mu|_{H_j(p)}$ denote the restriction of the measure $\mu$ to the hyperplane $H_j(p)$. Then, it holds that
\begin{equation}\label{eq:measure_extracted_form}
\mu|_{H_j(p)} = c_j(p)\pi\lambda_{\R^{n-1}},
\end{equation}
where the constant $c_j(p) \geq 0$ is given by the limit in \eqref{eq:c_constant_shifted} and $\lambda_{\R^{n-1}}$ denotes the Lebesgue measure on $\R^{n-1}$.
\end{thm}

\proof
Without loss of generality, we assume that $j = 1$, \ie we are considering the restriction of the measure $\mu$ to a hyperplane that lies orthogonal to the first coordinate axis. Together with this measure, we consider the Herglotz-Nevanlinna function $q$ whose triple of representing parameters is equal to $(0,\vec{0},\mu)$.

Let $\sigma := \mu|_{H_1(p)}$ be viewed as a Borel measure on $\R^{n-1}$, \ie
$$\diff\sigma(t_2,\ldots,t_n) = \diff\mu(p,t_2,\ldots,t_n).$$
We use the notation $\mu = \til{\mu} + \sigma$ and, by Theorem \ref{thm:intRep_Nvar}, we have that
\begin{multline}\label{eq:intRep_measure_split}
q(\vec{z}) = \frac{1}{\pi^n}\int_{\R^n}K_n(\vec{z},\vec{t})\diff\mu(\vec{t}) \\
= \frac{1}{\pi^n}\int_{\R^n}K_n(\vec{z},\vec{t})\diff\til{\mu}(\vec{t}) + \frac{1}{\pi^n}\int_{\R^{n-1}}K_n(\vec{z},(p,t_2,\ldots,t_n))\diff\sigma(t_2,\ldots,t_n).
\end{multline}

Investigating now the integral with respect to the measure $\sigma$, we calculate that
$$\begin{array}{LCL}
\multicolumn{3}{L}{\frac{1}{\pi^n}\int_{\R^{n-1}}K_n(\vec{z},(p,t_2,\ldots,t_n))\diff\sigma(p,t_2,\ldots,t_n)} \\[0.5cm]
~~~~~~~ & = & \frac{\I}{\pi^n}\int_{\R^{n-1}}\left(\frac{2}{(2\I)^n}\left(\frac{1}{p-z_1}-\frac{1}{p+\I}\right)\prod_{\ell = 2}^n\left(\frac{1}{t_\ell-z_\ell}-\frac{1}{t_\ell+\I}\right)\right. \\[0.5cm]
~ & ~ & \left.-\frac{1}{(2\I)^n}\left(\frac{1}{p-\I}-\frac{1}{p+\I}\right)\prod_{\ell = 2}^n\left(\frac{1}{t_\ell-\I}-\frac{1}{t_\ell+\I}\right)\right)\diff\sigma(t_2,\ldots,t_n) \\[0.5cm]
~ & = & \frac{\I}{\pi^n}\int_{\R^{n-1}}\left(\frac{2}{(2\I)^{n-1}}\:\frac{z_1+\I}{2\I(p-z_1)(p+\I)}\prod_{\ell = 2}^n\left(\frac{1}{t_\ell-z_\ell}-\frac{1}{t_\ell+\I}\right)\right. \\[0.5cm]
~ & ~ & \left.-\frac{1}{(2\I)^{n-1}}\:\frac{1}{1+p^2}\prod_{\ell = 2}^n\left(\frac{1}{t_\ell-\I}-\frac{1}{t_\ell+\I}\right)\right)\diff\sigma(t_2,\ldots,t_n) = (*).
\end{array}$$
The trick now is to add and subtract a term, such that the remaining expressions containing the $(t_2,\ldots,t_n)$-variables can be seen as the kernel $K_{n-1}$, multiplied by a factor independent of these variables. This leads to
$$\begin{array}{LCL}
(*) & = & \frac{\I}{\pi^n}\int_{\R^{n-1}}\left(\frac{2}{(2\I)^{n-1}}\:\frac{z_1+\I}{2\I(p-z_1)(p+\I)}\prod_{\ell = 2}^n\left(\frac{1}{t_\ell-z_\ell}-\frac{1}{t_\ell+\I}\right)\right. \\[0.5cm]
~ & ~ & \left.-\frac{1}{(2\I)^{n-1}}\:\frac{z_1+\I}{2\I(p-z_1)(p+\I)}\prod_{\substack{\ell = 1 \\ \ell \neq j}}^n\left(\frac{1}{t_\ell-\I}-\frac{1}{t_\ell+\I}\right)\right. \\[0.5cm]
~ & ~ & \left.+\frac{1}{(2\I)^{n-1}}\left(\frac{z_1+\I}{2\I(p-z_1)(p+\I)}-\frac{1}{1+p^2}\right)\right. \\[0.5cm]
~ & ~ & \cdot\left.\prod_{\ell = 2}^n\left(\frac{1}{t_\ell-\I}-\frac{1}{t_\ell+\I}\right)\right)\diff\sigma(t_2,\ldots,t_n) \\[0.5cm]
~ & = & \frac{1}{\pi}\left(\frac{z_1+\I}{2\I(p-z_1)(p+\I)}\cdot q_{1}(z_2,\ldots,z_n) + \frac{z_1-\I}{2\I(p-z_1)(p-\I)}\cdot q_{1}(\I,\ldots,\I)\right)
\end{array}$$
with the auxiliary function {$q_1$ being defined as}
\begin{multline}\label{eq:auxiliary}
q_{1}(z_2,\ldots,z_n) := \frac{\I}{\pi^{n-1}}\int_{\R^{n-1}}\left(\frac{2}{(2\I)^{n-1}}\prod_{\ell = 2}^n\left(\frac{1}{t_\ell-z_\ell}-\frac{1}{t_\ell+\I}\right)\right. \\
\left.-\frac{1}{(2\I)^{n-1}}\prod_{\ell = 2}^n\left(\frac{1}{t_\ell-\I}-\frac{1}{t_\ell+\I}\right)\right)\diff\sigma(t_2,\ldots,t_n) \\
= \frac{1}{\pi^{n-1}}\int_{\R^{n-1}}K_{n-1}((z_2,\ldots,z_n),(t_2,\ldots,z_n))\diff\sigma(t_2,\ldots,t_n).
\end{multline}

Note that, at this point, we do not know if the function $q_1$ is  well-defined. To this end, we observe now that
\begin{multline}\label{eq:measure_splitting_growth}
\int_{\R^n}\prod_{\ell=1}^n\frac{1}{1+t_\ell^2}\diff\mu(\vec{t}) \\
= \int_{\R^n}\prod_{\ell=1}^n\frac{1}{1+t_\ell^2}\diff\til{\mu}(\vec{t}) + \frac{1}{1+p^2}\int_{\R^{n-1}}\prod_{\ell = 2}^n\frac{1}{1+t_\ell^2}\diff\sigma(t_2,\ldots,t_n).
\end{multline}
All of the three terms above are non-negative since all of the integrands and  measures are positive. But the term on the left is finite, so both terms on the right must also be finite as well. This implies, first, that the function $q_{1}$ is well-defined and holomorphic in the poly-upper half-plane $\C^{+(n-1)}$ of dimension $n-1$. Furthermore, we are able to change the order of limits and integrations whenever we have either of the two measures on the right-hand side of equality \eqref{eq:measure_splitting_growth}, as we are allowed to do this with the measure on the left-hand side of equality \eqref{eq:measure_splitting_growth}.

Using this, we calculate that
$$\begin{array}{LCL}
\multicolumn{3}{L}{\lim\limits_{z_1 \ntto p}(p-z_1)\frac{1}{\pi^n}\int_{\R^n}K_n(\vec{z},\vec{t})\diff\til{\mu}(\vec{t}) = \frac{1}{\pi^n}\int_{\R^n} \lim\limits_{z_1 \ntto p}(p-z_1)K_n(\vec{z},\vec{t})\diff\til{\mu}(\vec{t})} \\[0.5cm]
~ & = & \frac{1}{\pi^n}\int_{\R^n}\frac{1}{(2\I)^{n-1}}\prod_{\ell = 2}^n\left(\frac{1}{t_\ell-z_\ell}-\frac{1}{t_\ell+\I}\right)\chi_{\{p\}}(t_1)\diff\til{\mu}(\vec{t}) \\[0.5cm]
~ & = & \frac{1}{\pi^n}\cdot\frac{1}{(2\I)^{n-1}}\int_{\R^{n-1}}\prod_{\ell = 2}^n\left(\frac{1}{t_\ell-z_\ell}-\frac{1}{t_\ell+\I}\right)\diff\til{\mu}(p,t_2,\ldots,t_n) = 0
\end{array}$$
as the measure $\til{\mu}$ is, by construction, identically zero on the hyperplane $H_j(p)$.

By Corollary \ref{coro:propeties_of_HN_funcitons}, there exists a number {$c := c_1(p) \geq 0$}, independent of the value of $(z_2,\ldots,z_n) \in \C^{+(n-1)}$, such that
$$\lim\limits_{z_1 \ntto p}(p-z_1)q(\vec{z}) = c.$$
But this limit is, by the expansion \eqref{eq:intRep_measure_split}, also equal to
$$\begin{array}{LCL}
\multicolumn{3}{L}{\lim\limits_{z_1 \ntto p}(p-z_1)q(\vec{z}) = \lim\limits_{z_1 \ntto p}(p-z_1)\frac{1}{\pi^n}\int_{\R^n}K_n(\vec{z},\vec{t})\diff\til{\mu}(\vec{t})} \\[0.5cm]
~ & ~ & + \lim\limits_{z_1 \ntto p}(p-z_1)\frac{1}{\pi}\left(\frac{z_1+\I}{2\I(p-z_1)(p+\I)}\cdot q_{1}(z_2,\ldots,z_n)\right. \\[0.5cm]
~ & ~ & +\left.\frac{z_1-\I}{2\I(p-z_1)(p-\I)}\cdot q_{1}(\I,\ldots,\I)\right) \\[0.5cm]
~& = & \frac{1}{2\pi\I}\left(q_{1}(z_2\ldots,z_n)+q_{1}(\I,\ldots,\I)\right).
\end{array}$$
We thus have that
$$c = \frac{1}{2\pi\I}\left(q_{1}(z_2\ldots,z_n)+q_{1}(\I,\ldots,\I)\right).$$
Since this {equality} holds for every vector $(z_2\ldots,z_n) \in \C^{+(n-1)}$, we can set {first} $(z_2\ldots,z_n) = (\I,\ldots,\I)$. This gives that $q_{1}(\I,\ldots,\I) = \pi\:\I\:c$, allowing us to solve the above equation for $q_{1}(z_2,\ldots,z_n)$, yielding
$$q_{1}(z_2,\ldots,z_n) = \pi\:\I\:c.$$

We infer now that the function $q_1$ is a Herglotz-Nevanlinna function in $n-1$ variables. On one side, its representing measure is $\mu|_{H_1(p)}$ due to equality \eqref{eq:auxiliary} and \cite[Corollar 4.7]{LugerNedic2017}. On the other hand, we know that the representing measure of the function $\vec{z} \to \I$, as a function of $n-1$ variables, is $\lambda_{\R^{n-1}}$, \cf \cite[Example 3.5]{LugerNedic2017}. Invoking the uniqueness statement of Theorem \ref{thm:intRep_Nvar} finishes the proof.
\endproof

The following corollary of Theorem \ref{thm:affine_hyperplane_extraction} gives a particular decomposition of {a} Herglotz-Nevanlinna function {with respect to a collection of coordinate-orthogonal hyperplanes}.

\begin{coro}\label{coro:affine_hyperplane_extraction}
Let $n \geq 2$, let $q$ be a Herglotz-Nevanlinna function and let $\mu$ be its representing measure. Decompose the measure $\mu$ as
$${\mu} =  \sum_{i\in I}\mu|_{H_{j_i}(p_i)}+\til{\mu}$$
for some indices $j_i \in \{1,2,\ldots,n\}$ and some points $p_i \in \R$, where $I {\subseteq \{1,2,\ldots,n\}}$ is a finite set of indices and $\til{\mu}$ denotes the remaining positive Borel measure. Then, the function $q$ can be written, for any $\vec{z} \in \C^{+n}$, as
\begin{equation}\label{eq:intRep_measure_one_extraction}
q(\vec{z}) = \sum_{i\in I}\frac{c_{j_i}(p_i)}{p_i-z_{j_i}} + \til{q}(\vec{z}),
\end{equation}
where $\til{q}$ is a Herglotz-Nevanlinna function which admits an integral representation formula of the form
\begin{equation}\label{eq:intRep_measure_one_extraction_part}
\til{q}(\vec{z}) = \left(a-\sum_{i\in I}\frac{{c_{j_i}(p_i)}\:p_i}{1+p_i^2}\right) + \sum_{\ell = 1}^{n}b_\ell z_\ell + \frac{1}{\pi^n}\int_{\R^n}K_n(\vec{z},\vec{t})\diff\til{\mu}(\vec{t})
\end{equation}
where $K_n$, $a$, $\vec{b}$ and $\mu$ are as in Theorem \ref{thm:intRep_Nvar} and $c_{j_i}(p_i)$ are given by the limit \eqref{eq:c_constant_shifted}.
\end{coro}

\begin{remark}
In the above corollary, we assume that the pairs $(j_i,p_i)$ are distinct in the sense that there do not exist $i_1,i_2 \in I$, such that both $j_{i_1} = j_{i_2}$ and $p_{i_1} = p_{i_2}$.
\end{remark}

\proof
Without loss of generality, we may assume that $I = \{1\}$ and write $p := p_1$. If this is the case, we {infer from} the proof of Theorem \ref{thm:affine_hyperplane_extraction} that the function $q$ can be written as
\begin{multline*}
q(z) =  a + \sum_{\ell=1}^n b_\ell\:z_\ell + \frac{1}{\pi^n}\int_{\R^n}K_n(\vec{z},\vec{t})\diff\til{\mu}(\vec{t}) \\
+ \frac{1}{\pi}\left(\frac{z_1+\I}{2\I(p-z_1)(p+\I)}\cdot q_{1}(z_2,\ldots,z_n) + \frac{z_1-\I}{2\I(p-z_1)(p-\I)}\cdot q_{1}(\I,\ldots,\I)\right),
\end{multline*}
where $q_1$ is the auxiliary function defined by formula {\eqref{eq:auxiliary}}. Furthermore, we have also learned that the function $q_1$ is, in fact, identically equal to $\pi\:\I\:c_1(p)$, yielding that
\begin{multline*}
 \frac{1}{\pi}\left(\frac{z_1+\I}{2\I(p-z_1)(p+\I)}\cdot q_{1}(z_2,\ldots,z_n) + \frac{z_1-\I}{2\I(p-z_1)(p-\I)}\cdot q_{1}(\I,\ldots,\I)\right) \\[0.35cm]
= \frac{c_1(p)\:(1+p\:z_1)}{(p-z_1)(1+p^2)} = \frac{c_1(p)}{p-z_1} - \frac{c_1(p)\:p}{1+p^2}.
\end{multline*}
This implies that expansion \eqref{eq:intRep_measure_split} takes the desired form \eqref{eq:intRep_measure_one_extraction}, with the function $\til{q}$ given, indeed, by representation \eqref{eq:intRep_measure_one_extraction_part}.

Thus, it remains to conclude that the function $\til{q}$, as defined in the theorem, is a Herglotz-Nevanlinna function. To do this, we only need to check that the measure $\til{\mu}$ satisfies the Nevanlinna condition \eqref{eq:measure_Nevan}, since we have already shown in the proof of Theorem \ref{thm:affine_hyperplane_extraction} that it satisfies the growth condition \eqref{eq:measure_growth} when we considered equality \eqref{eq:measure_splitting_growth}.

Without loss of generality, we restrict ourselves to only check the case $\ell_1 = 1$ and $\ell_2 = 2$ in condition \eqref{eq:measure_Nevan}. Hence, we consider the identity
\begin{multline}
\label{eq:measure_Nevan_splitting}
\int_{\R^n}\frac{1}{(t_1 - z_1)^2(t_2 - \bar{z_2})^2}\prod_{j=3}^n\left(\frac{1}{t_j - z_j} - \frac{1}{t_j - \bar{z}_j}\right)\diff\mu(\vec{t}) \\
= \int_{\R^n}\frac{1}{(t_1 - z_1)^2(t_2 - \bar{z_2})^2}\prod_{j=3}^n\left(\frac{1}{t_j - z_j} - \frac{1}{t_j - \bar{z}_j}\right)\diff\til{\mu}(\vec{t}) \\
+ \frac{C}{(p-z_1)^2}\cdot\int_{\R^{n-1}}\frac{1}{(t_2 - \bar{z_2})^2}\prod_{j=3}^n\left(\frac{1}{t_j - z_j} - \frac{1}{t_j - \bar{z}_j}\right)\diff t_2 \ldots \diff t_n,
\end{multline}
where $C \geq 0$ is some constant. The left-hand side of equality \eqref{eq:measure_Nevan_splitting} is identically zero for any vector $\vec{z} \in \C^{+n}$ since $\mu$ is the representing measure of a Herglotz-Nevanlinna function and, thus, satisfies the Nevanlinna condition \eqref{eq:measure_Nevan}. Furthermore, the second term on the right-hand side of equality \eqref{eq:measure_Nevan_splitting} is, likewise, identically zero since
$$\int_\R\frac{1}{(t_2-\bar{z_2})^2}\diff t_2 = 0$$
for any $z_2 \in \C^+$ by standard residue calculus. Therefore, the first term on the right-hand side of equality \eqref{eq:measure_Nevan_splitting} is also identically zero for any vector $\vec{z} \in \C^{+n}$, implying that the measure $\til{\mu}$ does indeed satisfy the Nevanlinna condition \eqref{eq:measure_Nevan} and finishing the proof.
\endproof

Note that in one variable the decomposition of the function in formula \eqref{eq:intRep_measure_one_extraction} in Corollary \ref{coro:affine_hyperplane_extraction} is, of course, well known. Given just a single point $p \in \R$, one can find a number $c \geq 0$, such that $$\mu|_{\{p\}} = c\pi\delta_p,$$
which yields a decomposition of the function $q$ as
$$q(z) = \frac{c}{p-z} + \til{q}(z),$$
where $\til{q}$ is a Herglotz-Nevanlinna function whose representing measure is $\mu-\mu|_{\{p\}}$. Note that the convention of writing the measure $\mu|_{\{p\}}$ as $c\pi\delta_p$, and not $c'\delta_p$, comes from the fact that the Herglotz-Nevanlinna function $z \mapsto -\frac{1}{z}$ is represented by the measure $\pi\delta_0$.

In several variables, the procedure of decomposing a Nevanlinna measure and obtaining a decomposition of a Herglotz-Nevanlinna function as in formula \eqref{eq:intRep_measure_one_extraction} cannot be generalized to arbitrary sets. Indeed, given a Nevanlinna measure $\mu$ and given any Borel measurable subset $U \subseteq\R^n$, one might consider the decomposition
\begin{equation}\label{eq:measure_decomposition}
\mu = \mu|_U + (\mu-\mu|_U).
\end{equation}
In one variable, the only constraint on a representing measure is the growth condition \eqref{eq:measure_growth}, and, hence, both measures on the right-hand side of  equality \eqref{eq:measure_decomposition} are automatically Nevanlinna measures. However, in several variables this is not necessarily the case as the Nevanlinna condition \eqref{eq:measure_Nevan} must be fulfilled as well.

Let us now return to the $\mu$-mass of certain subsets of $\R^n$. Suppose so that $U \subseteq \R^n$ is a Borel measurable set, such that $U \subseteq H_j(p)$ for some index $j \in \{1,2\ldots,n\}$ and some point $p \in \R$. Then, obviously,
$$\mu(U) = 0 \iff \mu|_{H_j(p)}(U) = 0.$$
A direct implication of this trivial fact, together with Theorem \ref{thm:affine_hyperplane_extraction}, is the following important statement.

\begin{coro}\label{coro:measure_affine_subspaces}
Let $n \geq 2$ and let $\mu$ be a Nevanlinna measure on $\R^n$. Let $U$ be an affine subspace of $\R^n$ which is orthogonal to some coordinate axis. Then it holds: if $\codim(U) \geq 2$ then $\mu(U)$ is zero, {while} if $\codim(U) =1 $ then {$\mu(U)$} is either zero or infinity.
\end{coro}

In particular, for $n \geq 2$, points, \ie affine subspaces of $\codim(U) =n$, have measure zero.

From the above corollary, we conclude that if an affine subspace which is orthogonal to some coordinate axis is to hope to have non-zero $\mu$-mass, it needs to have codimension one. Even then, the only non-zero option is infinity, and the measure on the subspace may only have the form of the Lebesgue measure due to Theorem \ref{thm:affine_hyperplane_extraction}.

The final corollary of this section establishes a relation between the variable dependence of the function $q$ and the $\mu$-masses of the hyperplanes $H_j(0)$.

\begin{coro}\label{coro:variable_dependence}
Let $n \geq 2$, let $q$ be a Herglotz-Nevanlinna function and let $\mu$ be its representing measure. If $\mu(H_j(0)) = \infty$, then $q$ has to depend on the $z_j$-variable.
\end{coro}

\proof
As in the proof of Theorem \ref{thm:affine_hyperplane_extraction}, we may, without loss of generality, assume that $j = 1$. Let us now do a proof by contradiction. Suppose so that $\mu(H_1(0)) = \infty$ and that the function $q$ does not depend on the $z_1$-variable.

By Theorem \ref{thm:affine_hyperplane_extraction}, we have that $\mu|_{H_1(0)} = c_1(0)\pi\lambda_{\R^{n-1}}$, where $c_1(0) \neq 0$ due to our assumption. On the other hand, Corollary \ref{coro:propeties_of_HN_funcitons} implies that
$$c_1(0) = -\lim\limits_{z_1 \ntto 0}z_1\:q(\vec{z}) = \lim\limits_{z_1 \ntto 0}z_1\:q(\I,z_2,\ldots,{z_n}) = 0.$$
This gives the desired contradiction, finishing the proof. 
\endproof

Observe, though, that the converse to Corollary \ref{coro:variable_dependence}, \ie that functions depending on the $z_j$-variable have to have $\mu(H_j(0)) = \infty$, is not true, as demonstrated by the function $q(z_1,z_2) = -\tfrac{1}{z_1+z_2}$.

\subsection{Geometry of the support}\label{subsec:support}

As we have seen before, hyperplanes which are orthogonal to some coordinate axis appear as support sets for Nevanlinna measures. However, for rotated hyperplanes, this may or may not be true.

\begin{example}\label{ex:lines}
Let $n = 2$ and consider the following three hyperplanes in $\R^2$: $H_1 := \{t_1 = 1\}$, $H_2 := \{t_1 = -t_2\}$ and $H_3 := \{t_1 = t_2\}$.

We know from Theorem \ref{thm:affine_hyperplane_extraction} that the hyperplane $H_1$ appears as the support of some Nevanlinna measure. Similarly, for the hyperplane $H_2$, one can find a Nevanlinna measure whose support is equal to $H_2$, \cf Example \ref{ex:lines_v2} and \cite[Example 4.2]{Nedic2017}. However, we will soon see that there exists no Nevanlinna measure such that its support would be equal to, {or even contained in,} $H_3$, \cf Theorem \ref{thm:non_par_hyperplanes}, Example \ref{ex:diagonal} and Example \ref{ex:lines_v2}. \hfill$\lozenge$
\end{example}

In what follows, we are interested in identifying subsets  of $\R^n$ which cannot contain the support of a Nevanlinna measure.

\begin{thm}
\label{thm:non_par_hyperplanes}
Let $n \geq 2$ and let $\mu$ be a Nevanlinna measure. Define the affine subspace $\Sigma(\mat{A},\vec{\beta}) \subseteq \R^n$ as
$$\Sigma(\mat{A},\vec{\beta}) := \{\mat{A}\vec{s}+\vec{\beta}~|~\vec{s} \in \R^n\},$$
where $\mat{A}  = \{\alpha_{i,j}\}_{i,j=1}^n \in M_n(\R)$ and $\vec{\beta} \in \R^n$ are such that the matrix $\mat{A}$ {contains no trivial} rows and that there exist {two distinct} indices $j_1,j_2 \in \{1,2,\ldots,n\}$ and a number $\gamma > 0$, such that $\alpha_{j_1,\ell} = \gamma\:\alpha_{j_2,\ell}$ for all $\ell \in \{1,2,\ldots,n\}$.
Then,
\begin{equation}
\label{eq:non_par_hyperplanes}
\supp(\mu) \subseteq \Sigma(\mat{A},\vec{\beta}) \Longrightarrow \mu \equiv 0.
\end{equation}
\end{thm}

\begin{remark}
The conditions on the matrix $A$ say, in other words, that the set $\Sigma(\mat{A},\vec{\beta})$, firstly, should not be coordinate-parallel and, secondly, it should be contained in a hyperplane of the form $\gamma x_{j_1}=x_{j_2}$ for some indices $j_1$ and $j_2$ and a positive dependence-factor $\gamma$. In the case $n=2$, these sets are precisely lines with positive slope.
\end{remark}

\proof
Let $\supp(\mu) \subseteq \Sigma(\mat{A},\vec{\beta})$ for $\mat{A}$ and $\vec{\beta}$ as in the theorem. Without loss of generality, we may assume that $j_1 = 1$ and $j_2 = 2$. Since $\mu$ is a Nevanlinna measure, it satisfies, in particular, the condition that
\begin{equation}
\label{eq:nevan_condition_v2}
\int_{\R^n}\frac{1}{(t_1 - z_1)^2(t_2 - \bar{z_2})^2}\prod_{j=3}^n\left(\frac{1}{t_j - z_j} - \frac{1}{t_j - \bar{z}_j}\right)\diff\mu(\vec{t}) = 0
\end{equation}
for any $\vec{z} \in \C^{+n}$. In our case, this integral may be rewritten as
\begin{multline}
\label{eq:nevan_int_1}
\int_{\R^n}\frac{1}{((\mat{A}\vec{s}+\vec{\beta})_1 - z_1)^2((\mat{A}\vec{s}+\vec{\beta})_2 - \bar{z_2})^2} \\
\cdot\prod_{j=3}^n\left(\frac{1}{(\mat{A}\vec{s}+\vec{\beta})_j - z_j} - \frac{1}{(\mat{A}\vec{s}+\vec{\beta})_j - \bar{z_j}}\right)\diff\mu(\mat{A}\vec{s}+\vec{\beta}) = 0.
\end{multline}

The assumptions on the matrix $\mat{A}$ give now that
$$(\mat{A}\vec{s})_2 = \gamma\:(\mat{A}\vec{s})_1,$$
yielding further that
$$\frac{1}{((\mat{A}\vec{s}+\vec{\beta})_1 - z_1)^2((\mat{A}\vec{s}+\vec{\beta})_2 - \bar{z_2})^2} = \frac{1}{\gamma^2((\mat{A}\vec{s})_1 + \beta_1 - z_1)^2((\mat{A}\vec{s})_1 + \frac{\beta_2}{\gamma}- \frac{\bar{z_2}}{\gamma})^2}.$$
Choosing now the point
$$(\I,\beta_2-\beta_1\gamma+\gamma\:\I,\I,\ldots,\I) \in \C^{+n},$$
we calculate that the integral \eqref{eq:nevan_int_1} at this point is equal to
$$(2\I)^{n-2}\int_{\R^n}\frac{1}{\gamma^2|(A\vec{s})_1+\beta_1-\I|^4}\prod_{j=3}^n\frac{1}{1+((\mat{A}\vec{s})j+\beta_j)^2}\diff\mu(\mat{A}\vec{s}+\vec{\beta}) = 0,$$
{and} since the integrand is a positive function and $\supp(\mu) \subseteq \Sigma(\mat{A},\vec{\beta})$, we must have $\mu \equiv 0$. This finishes the proof.
\endproof

\begin{remark}
If the matrix $\mat{A}$ in the formulation of Theorem \ref{thm:non_par_hyperplanes} would contain a trivial row, then the set $\Sigma(\mat{A},\vec{\beta})$ would be contained in a hyperplane orthogonal to some coordinate axis. Thus, it is covered by Theorem \ref{thm:affine_hyperplane_extraction} and Corollary \ref{coro:measure_affine_subspaces}.
\end{remark}

Let us now consider some examples which show the use of Theorem \ref{thm:non_par_hyperplanes} as well as the necessity of its requirements on the matrix $\mat{A}$.

\begin{example}\label{ex:diagonal}
Let us choose
$$\mat{A} = \left[\begin{array}{cccc}
1 & 0 & \dots & 0 \\
1 & 0 & \dots & 0 \\
\vdots & \vdots & \ddots & \vdots \\
1 & 0 & \dots & 0
\end{array}\right]_{n \times n}.$$
Then, the set $\Sigma(\mat{A},\vec{0})$ equals the diagonal in $\R^n$, with Theorem \ref{thm:non_par_hyperplanes} showing that a Nevanlinna measure cannot be supported only in this set.\hfill$\lozenge$
\end{example}

\begin{example}
Let $n = 3$ and let us choose
$$\mat{A} = \left[\begin{array}{ccc}
1 & 0 & 0 \\
0 & 1 & 0 \\
-1 & -1 & 0
\end{array}\right].$$
Then, the set $\Sigma(\mat{A},\vec{0})$ equals the plane $\{t_1 + t_2 + t_3 = 0\} \subseteq \R^3$. However, the matrix $A$, though not of maximal rank, does not satisfy the assumption of Theorem \ref{thm:non_par_hyperplanes} about having a pair of linearly dependent rows. Thus, Theorem \ref{thm:non_par_hyperplanes} does not apply, and, in fact, it can be shown that the set $\Sigma(\mat{A},\vec{0})$ equals the support of the representing measure of the Herglotz-Nevanlinna function
$$(z_1,z_2,z_3) \mapsto \frac{-1}{z_1+z_2+z_3},$$
see also \cite[Example 4.7]{LugerNedic2017}. \hfill$\lozenge$
\end{example}

\begin{example}\label{ex:lines_v2}
Let $n = 2$ and let us choose
$$\mat{A} = \left[\begin{array}{cc}
1 & 0 \\
\xi & 0
\end{array}\right] \quad\text{and}\quad \vec{\beta} = \left[\begin{array}{c}
0 \\
\eta
\end{array}\right].$$
Then, the set $\Sigma(\mat{A},\vec{\beta})$ equals the line $\{(t_1,t_2) \in \R^2~|~\xi t_1 + \eta = t_2\} \subseteq \R^2$. If $\xi > 0$, then, by Theorem \ref{thm:non_par_hyperplanes}, the set $\Sigma(\mat{A},\vec{\beta})$ cannot contain the support of some Nevanlinna measure. However, if $\xi < 0$, Theorem \ref{thm:non_par_hyperplanes} does not apply, and one can, in fact, find a Herglotz-Nevanlinna fucntion such that the support of its representing measure is equal to $\Sigma(\mat{A},\vec{\beta})$. In particular, when $\xi = -1$ and $\eta = 0$, the set $\Sigma(\mat{A},\vec{\beta})$ equals the anti-diagonal in $\R^2$, which equals the support of the representing measure of the Herglotz-Nevanlinna function
$$(z_1,z_2) \mapsto \frac{-1}{z_1+z_2},$$
see also \cite[Example 4.2]{Nedic2017}.\hfill$\lozenge$
\end{example}

\begin{figure}[!ht]
\includegraphics[width=0.5\linewidth]{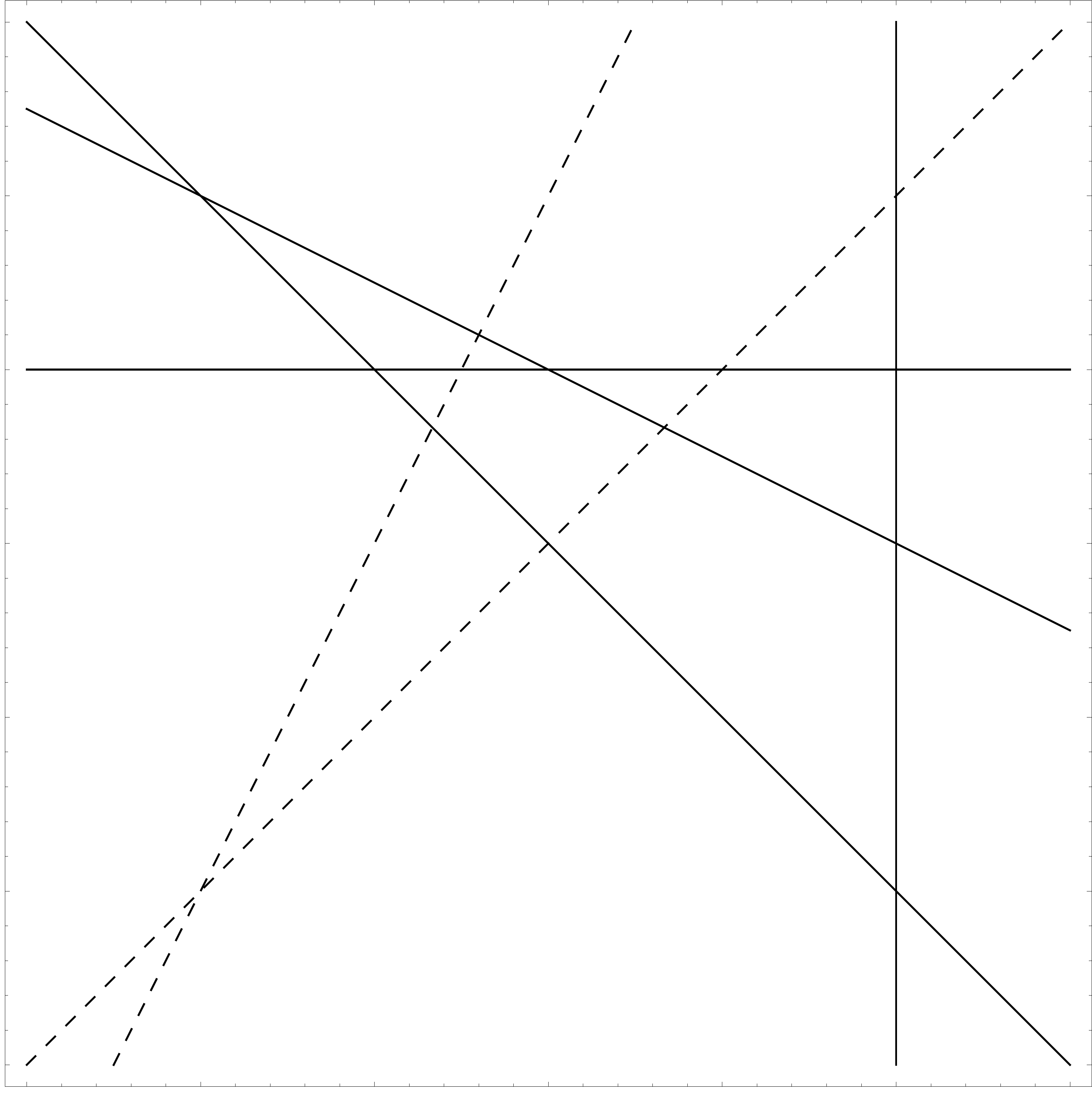}
\caption{The solid lines are those which {can contain} the support of some Nevanlinna measure, while the dashed lines are those {which} cannot contain the support of any Nevanlinna measure, \cf Examples \ref{ex:lines} and \ref{ex:lines_v2}. The plot area is $[-3,3]^2 \subseteq \R^2$.}\label{fig:lines}
\end{figure}

More generally, we show now that the support of a Nevanlinna measure cannot {even} be confined  to a strip of positive slope.

\begin{thm}
\label{thm:strips}
Let $n \geq 2$ and let $\mu$ be a Nevanlinna measure. Define the strip
$$S_{j_1,j_2}(\alpha,\beta_1,\beta_2) := \{\vec{t}\in \R^n~|~\beta_1 < t_{j_2} - \alpha\:t_{j_1} < \beta_2 \},$$
where $\alpha > 0$ and $\beta_1,\beta_2 \in \R$ are two numbers with $\beta_1 < \beta_2$, and $j_1,j_2 \in \{1,2,\ldots,n\}$ two distinct indices.
Then,
$$\supp(\mu) \subseteq S_{j_1,j_2}(\alpha,\beta_1,\beta_2) \Longrightarrow \mu \equiv 0.$$
\end{thm}

\proof
Let $\alpha, \beta_1$ and $\beta_2$ be as in the theorem. Without loss of generality, we may assume that $j_1 = 1$ and $j_2 = 2$. We will now show that there exists a vector $\vec{z} \in \C^{+n}$, for which condition \eqref{eq:nevan_condition_v2} is fulfilled only in the case of the zero measure.

First, we observe that
\begin{multline*}
\frac{1}{(t_1-z_1)^2(t_2-\bar{z_2})^2} = \frac{(t_1-\bar{z_1})^2(t_2-z_2)^2}{|t_1-z_1|^4|t_2-\bar{z_2}|^4} \\
= \frac{((t_1-x_1)(t_2-x_2)+y_1y_2)^2-((t_1-x_1)y_2-(t_2-x_2)y_1)^2}{|t_1-z_1|^4|t_2-\bar{z_2}|^4} \\
-2\I\frac{((t_1-x_1)(t_2-x_2)+y_1y_2)((t_1-x_1)y_2-(t_2-x_2)y_1)}{|t_1-z_1|^4|t_2-\bar{z_2}|^4},
\end{multline*}
where we used $x_i := \Re[z_i], y_i := \Im[z_i]$ for $i = 1,2$. Introduce now two new variables $s_1$ and $s_2$ as $s_1 := \frac{1}{2}(\alpha\:t_1 + t_2)$ and $s_2 := \frac{1}{2}(\alpha\:t_1-t_2)$. Choose  $y_1, y_2 > 0$ such that $y_2 = \alpha\: y_1$ and denote $y := y_2 = \alpha\:y_1$. With these choices, the numerator of the imaginary part of the above expression becomes
\begin{multline*}
((t_1-x_1)(t_2-x_2)+y_1y_2)((t_1-x_1)y_2-(t_2-x_2)y_1) \\
= ((\tfrac{s_1+s_2}{\alpha}-x_1)(s_1-s_2-x_2)+\alpha\:y^2)((\tfrac{s_1+s_2}{\alpha}-x_1)y - \frac{1}{\alpha}(s_1-s_2-x_2)y) \\
= \tfrac{y}{\alpha^2}((s_1+s_2-\alpha\:x_1)(s_1-s_2-x_2)+\alpha^2y^2)(2s_2-\alpha\:x_1+x_2).
\end{multline*}
Since $\beta_1 < t_2-\alpha\:t_1 < \beta_2$, we infer that $\beta_1 < -2s_2 < \beta_2$, implying that the parameter $s_2$ is bounded. As such, we may choose $x_1$ and $x_2$ such that $2s_2-\alpha\:x_1+x_2 > 0$.

We may now adjust the choice of $y$ to ensure that the expression $(s_1+s_2-\alpha\:x_1)(s_1-s_2-x_2)+\alpha^2y^2$ is also positive. Indeed, observe first that
\begin{multline*}
(s_1+s_2-\alpha\:x_1)(s_1-s_2-x_2)+\alpha^2y^2 \\
= (s_1 - \tfrac{\alpha x_1 + x_2}{2})^2 - (\tfrac{\alpha x_1 + x_2}{2})^2 - (s_2-\alpha\:x_1)(s_2+x_2) + \alpha^2\:y^2.
\end{multline*}
The second and third term in the above expression depend only on $x_1,x_2$, which have been fixed, and $s_2$, which is bounded. Therefore, the value of $y$ may be chosen such that the sum is positive for all possible values of $s_2$.

In conclusion, we have shown that there exist values $x_1,x_2,y_1$ and $y_2$, such that
$$\frac{((t_1-x_1)(t_2-x_2)+y_1y_2)((t_1-x_1)y_2-(t_2-x_2)y_1)}{|t_1-z_1|^4|t_2-\bar{z_2}|^4} > 0$$
for all $t_1,t_2 \in S_{1,2}(\alpha,\beta_1,\beta_2)$. This means that, when considering the condition \eqref{eq:nevan_condition_v2} at the point
$$(x_1+\I\:y_1,x_2+\I\:y_2,\I,\ldots,\I) \in \C^{+n},$$
we have that the imaginary part of the term
$$\frac{1}{(t_1-z_1)^2(t_2-\bar{z_2})^2}$$
is always {negative}, while the term
$$\prod_{j=3}^n\left(\frac{1}{t_j - z_j} - \frac{1}{t_j - \bar{z}_j}\right) = (2\I)^{n-2}\prod_{j=3}^n\frac{1}{1+t_j^2}$$
is either pure-real or pure-imaginary. Depending on which option occurs here, we conclude that the integrand in condition \eqref{eq:nevan_condition_v2} either has non-zero real or imaginary part at this particular point in $\C^{+n}$, implying that the measure $\mu$ must be identically zero. This finishes the proof.
\endproof

\begin{example}
Let $n = 2$ and consider the strip
$$S_{1,2}(1,-1,1) = \{(t_1,t_2) \in \R^2~|~-1 < t_2-t_1 < 1\}.$$
Theorem \ref{thm:strips} now says that any non-zero Borel measure $\mu$ with support contained in this strip is not a Nevanlinna measure. In particular, we observe, again, that any non-zero measure {whose support is contained in the diagonal in $\R^2$} is not a Nevanlinna measure. \hfill$\lozenge$
\end{example}

\subsection{Refinements using M\"obius transforms}

Given a Herglotz-Nevanlinna function, one can precompose it with a product of M\"obius transforms that fix the upper half plane to obtain, again, a function of the same class. This can be used in order to obtain {a generalization of Corollary \ref{coro:affine_hyperplane_extraction}, as well as} further refinements of the geometric restrictions on the support of a Nevanlinna measure given by Theorems \ref{thm:non_par_hyperplanes} and \ref{thm:strips}. We give some such statements below. To start with we need some notation.

Consider a Nevanlinna measure $\mu$, such that for some point $\vec{p} \in \R^n$, it holds that $\mu|_{H_j(p_j)} \equiv 0$ for all $j = 1,2,\ldots,n$. Such a measure can be considered as a measure on {$\bigtimes_{j = 1}^n(\R\setminus\{p_j\})$} and, therefore, we may, for any collection of indices $\{j_1,j_2,\ldots,j_k\} \subseteq \{1,2,\ldots,n\}$, do $k$ changes of variables

\begin{equation}\label{eq:change_of_variables}
t_{j_\ell} \longmapsto \frac{1}{p_{j_\ell}-t_{j_\ell}}
\end{equation}

for $\ell = 1,2,\ldots,k$. Without loss of generality, we may restrict ourselves to investigate the case $k = 1$, $j_1 = 1$ and $p_1 = 0$. In this case, it suffices to assume that only $\mu|_{H_1(0)} \equiv 0$.

For a Borel set $U \subseteq (\R\setminus\{0\})\times\R^{n-1}$, we define
\begin{equation}\label{eq:map_J}
J_1^0(U) := \big\{(x_1,x_2,\ldots,x_n) \in \R^n~|~(-\tfrac{1}{x_1},x_2,\ldots,x_n) \in U\big\},
\end{equation}
{where the subscript $(\:\cdot\:)_1$ refers to taking $j_\ell = 1$ in formula \eqref{eq:change_of_variables}, while the superscript $(\:\cdot\:)^0$ refers to taking $p_{j_\ell} = 0$}. Similarly, for a Nevanlinna measure $\mu$ with $\mu|_{H_1(0)} \equiv 0$, we define
\begin{equation}\label{eq:map_J_star}
((J_1^0)^*\mu)(U) := \mu(J_1^0(U))
\end{equation}
for any Borel  set $U$ as before. The following proposition now justifies the introduction of these maps.

\begin{prop}\label{prop:measure_transform}
Let $\mu$ be a Nevanlinna measure with $\mu|_{H_1(0)} \equiv 0$.  Then, the measure $(J_1^0)^*\mu$ is also a Nevanlinna measure.
\end{prop}

\proof
We begin by investigating what happens to the integrals
$$\int_{\R\setminus\{0\}}\frac{1}{1+t^2}\diff t \quad\text{and}\quad \int_{\R\setminus\{0\}}\left(\frac{1}{t-z}-\frac{1}{t-w}\right)\diff t,$$
where $z,t \in \C\setminus\R$, under the change of variables $t = -\frac{1}{\tau}$. As such, we calculate that
$$\int_{\R\setminus\{0\}}\frac{1}{1+t^2}\diff t = \int_{\R\setminus\{0\}}\frac{1}{1+\tau^2}\diff \tau$$
and that
\begin{multline*}
\int_{\R\setminus\{0\}}\left(\frac{1}{t-z}-\frac{1}{t-w}\right)\diff t = \int_{\R\setminus\{0\}}\left(\frac{1}{-\frac{1}{\tau}-z}-\frac{1}{-\frac{1}{\tau}-w}\right)\frac{1}{\tau^2}\diff \tau \\
= \int_{\R\setminus\{0\}}\left(\frac{1}{\tau + \frac{1}{z}}-\frac{1}{\tau + \frac{1}{w}}\right)\diff \tau.
\end{multline*}

The first of the above calculations implies immediately that the measure $(J_1^0)^*\mu$ satisfies the growth condition \eqref{eq:measure_growth}, while the second calculation shows that, due to $z \mapsto -\frac{1}{z}$ being an automorphism of $\C^+$, the measure $(J_1^0)^*\mu$ also satisfies the Nevanlinna condition \eqref{eq:measure_Nevan}. The result then follows.
\endproof

Using {Proposition \ref{prop:measure_transform}}, the results of {Corollary \ref{coro:affine_hyperplane_extraction} and} Theorems \ref{thm:non_par_hyperplanes} and \ref{thm:strips} can be extended as follows.

\begin{coro}\label{coro:minus_representation}
{Let $n \geq 2$, let $q$ be a Herglotz-Nevanlinna function and let $\mu$ be its representing measure. Then, the Herglotz-Nevanlinna function
$$Q_1\colon(z_1,z_2,\ldots,z_n) \mapsto q(-\tfrac{1}{z_1},z_2,\ldots,z_n)$$
can be written, for any $\vec{z} \in \C^n$, as
$$Q_1(\vec{z}) = a + c_1(0)\:z_1 - \frac{b_1}{z_1} + \sum_{\ell = 2}^nb_\ell\:z_\ell + \frac{1}{\pi^n}\int_{\R^n}K_n(\vec{z},\vec{t})\diff\big((J_1^0)^*\til{\mu}\big)(\vec{t}),$$
where $K_n$, $a$, $\vec{b}$ and $\mu$ are as in Theorem \ref{thm:intRep_Nvar}, the number $c_{1}(0)$ is given by the limit \eqref{eq:c_constant_shifted} and
$$\til{\mu} :=  \mu - \mu|_{H_1(0)}.$$}
\end{coro}

\begin{coro}\label{coro:minus_inverse}
Let $n \geq 2$ and let $\mu$ be a Nevanlinna measure with $\mu|_{H_1(0)} \equiv 0$. Then, for sets $\Sigma(\mat{A},\vec{\beta})$ and $S_{j_1,j_2}(\alpha,\beta_1,\beta_2)$ as in Theorems \ref{thm:non_par_hyperplanes} and \ref{thm:strips}, {respectively,} it holds that
$$\supp(\mu) \subseteq J_1^0(\Sigma(\mat{A},\vec{\beta}) \cap (\R\setminus\{0\})\times\R^{n-1}) \Longrightarrow \mu \equiv 0$$
and
$$\supp(\mu) \subseteq J_1^0(S_{j_1,j_2}(\alpha,\beta_1,\beta_2) \cap (\R\setminus\{0\})\times\R^{n-1}) \Longrightarrow \mu \equiv 0.$$
\end{coro}

\begin{remark}
Any combination of maps $J_j^{p_j}$ and $(J_j^{p_j})^*$, defined analogously as the maps $J_1^0$ and $(J_1^0)^*$ in formulas \eqref{eq:map_J} and \eqref{eq:map_J_star}, respectively, can be used to extend the results of {Proposition \ref{prop:measure_transform} and} Corollar{ies \ref{coro:minus_representation} and} \ref{coro:minus_inverse}.

In particular, any set as in Theorems \ref{thm:non_par_hyperplanes} and \ref{thm:strips} may be successively transformed by maps of the from $J_j^{p_j}$, thereby enlarging or collection of subsets of $\R^n$ which cannot {contain} the support of some Nevanlinna measure, \cf Example \ref{ex:strips}.
\end{remark}

\begin{example}
Let $n = 2$ and consider the situation of Example \ref{ex:diagonal}, \ie choose
$$\mat{A} = \left[\begin{array}{cc}
1 & 0 \\
1 & 0
\end{array}\right].$$
Then, the set $\Sigma(\mat{A},\vec{0})$ equals the diagonal in $\R^2$. One can then calculate that
$$J_1^0(\Sigma(\mat{A},\vec{0}) \cap (\R\setminus\{0\})\times\R) = \{(\tau_1,\tau_2) \in \R^2~|~\tau_2 = -\tfrac{1}{\tau_1}\},$$
with Corollary \ref{coro:minus_inverse} now implying that a Nevanlinna measure $\mu$ cannot have its support contained only in the hyperbola given by the equation $\tau_2 = -\tfrac{1}{\tau_1}$. \hfill$\lozenge$
\end{example}

\begin{figure}[!ht]
\includegraphics[width=0.8\linewidth]{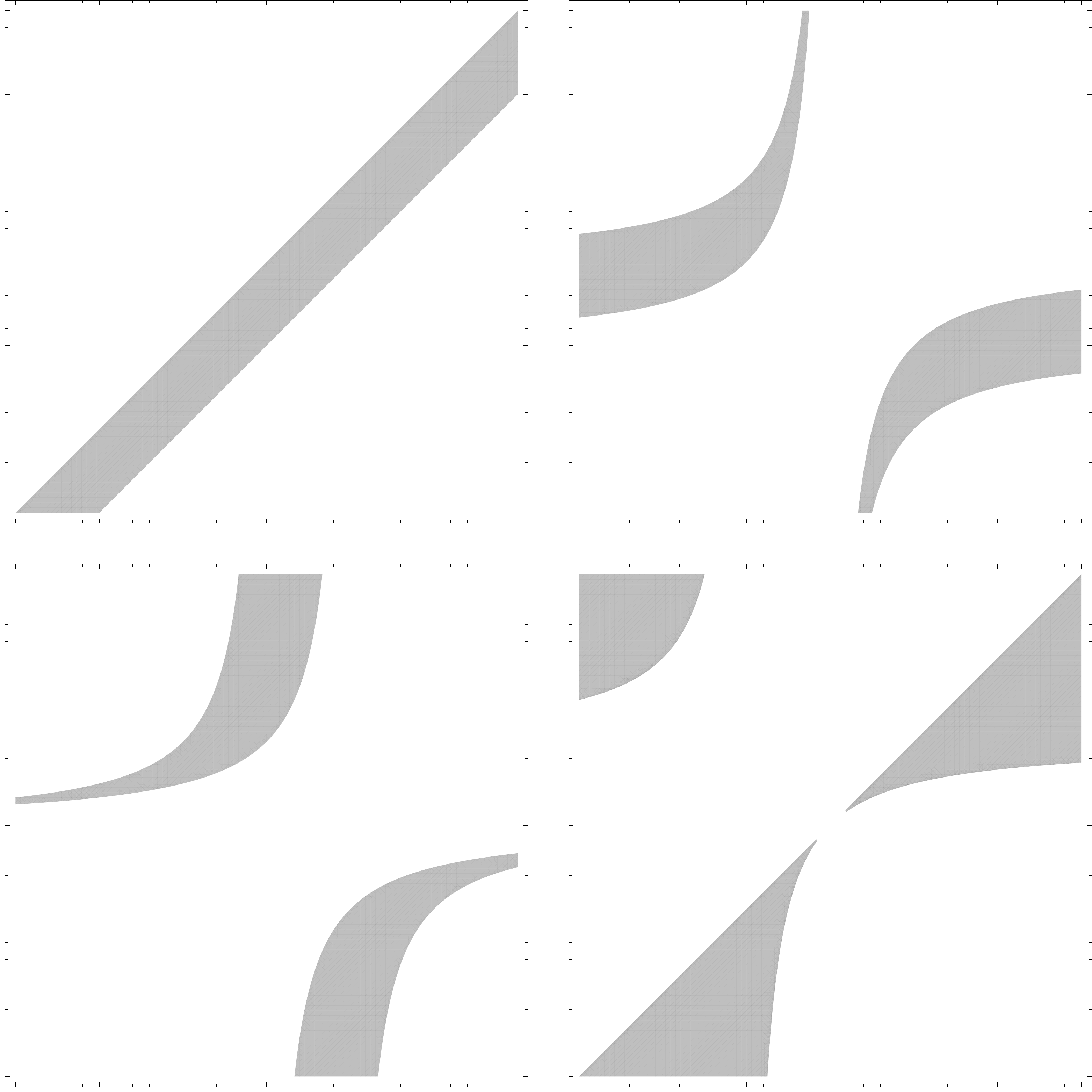}
\caption{The strip $S_{1,2}(1,-1,0)$ from Example \ref{ex:strips} (top left) and its transformations using the maps $J_1^0$ (top right), $J_2^0$ (bottom left) and $J_1^0J_2^0$ (bottom right). The plot area is always $[-3,3]^2 \subseteq \R^2$.}\label{fig:strips}
\end{figure}

\begin{example}\label{ex:strips}
Let $n = 2$ and consider the strip
$$S_{1,2}(1,-1,0) = \{(t_1,t_2) \in \R^2~|~0<t_1-t_2<1\}.$$
One can then establish its transformations using the maps $J_1^0$ and $J_2^0$ to be equal to
$$\begin{array}{rcl}
J_1^0(S_{1,2}(1,-1,0) \cap (\R\setminus\{0\})\times\R) & = & \{(\tau_1,\tau_2) \in \R^2~|~ 0 < -\tfrac{1}{\tau_1} - \tau_2 < 1\}, \\
J_2^0(S_{1,2}(1,-1,0) \cap (\R \times \R\setminus\{0\})) & = & \{(\tau_1,\tau_2) \in \R^2~|~ 0 < \tau_1 + \tfrac{1}{\tau_2} < 1\}, \\
J_1^0J_2^0(S_{1,2}(1,-1,0) \cap (\R\setminus\{0\})^2) & = & \{(\tau_1,\tau_2) \in \R^2~|~ 0 < -\tfrac{1}{\tau_1} + \tfrac{1}{\tau_2} < 1\}.
\end{array}$$
Theorem \ref{thm:strips} and Corollary \ref{coro:minus_inverse} imply now that a Nevanlinna measure $\mu$ cannot have its support contained in any of the above subsets of $\R^2$, \cf Figure \ref{fig:strips}. \hfill$\lozenge$
\end{example}

\begin{figure}[!ht]
\includegraphics[width=0.5\linewidth]{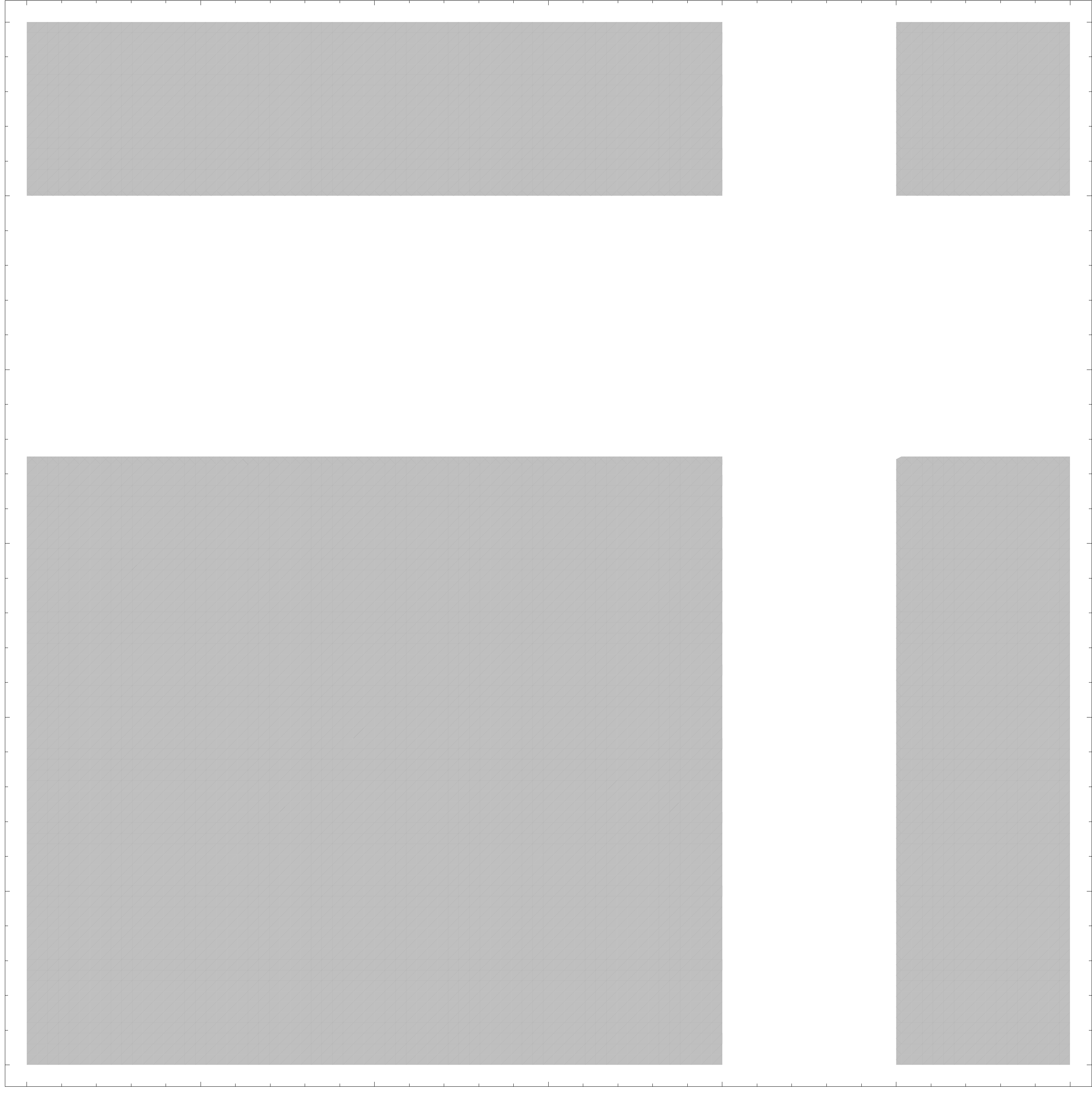}
\caption{The only Nevanlinna measure whose support {is contained} in the shown set is the trivial measure, \cf Theorem \ref{thm:cross}. This figure shows {the complement of the union of the} strips $\{(t_1,t_2) \in \R^2~|~1 < t_1 < 2\}$ and $\{(t_1,t_2) \in \R^2~|~\frac{1}{2} < t_2 < 2\}$. The plot area is $[-3,3]^2 \subseteq \R^2$.}\label{fig:cross}
\end{figure}

From Proposition \ref{prop:infinite_measure}, it is clear that the support of a Nevanlinna measure cannot be a bounded set. As a consequence of this fact and the  technique of coordinate transformation presented {previously in this section}, we can show that the support cannot be localized too much, in the sense that there cannot exist $n$ coordinate-orthogonal strips that do not intersect the support, \cf Figure \ref{fig:cross}.

\begin{thm}\label{thm:cross}
Let $n \geq 2$ and let $\mu$ be a Nevanlinna measure and suppose there exist numbers $\alpha_j < \beta_j$ for $j = 1,2,\ldots,n$, such that
$$\supp(\mu) \cap \bigcup_{j=1}^n\{\vec{t} \in \R^n~|~\alpha_j < t_j < \beta_j\} = \emptyset.$$
Then, $\mu \equiv 0$.
\end{thm}

\proof
Without loss of generality, we may assume that $n = 2$. Furthermore, we may also assume that $\alpha_1 = \alpha_2 = -1$ and $\beta_1 = \beta_2 = 1$ as other cases may be covered by scalings and translations.

As such, we suppose that
$$\supp(\mu) \cap \{(t_1,t_2) \in \R^2~|~-1 < t_1 < 1\} \cap  \{(t_1,t_2) \in \R^2~|~-1 < t_2 < 1\}= \emptyset.$$
In this case, the support of the measure $(J_1^0)^*(J_2^0)^*\mu$ is contained in the square $[-1,1]^2 \subseteq \R^2$, {and} is, therefore, a finite Nevanlinna measure. But, by Proposition \ref{prop:infinite_measure}, it now holds that $(J_1^0)^*(J_2^0)^*\mu \equiv 0$. This translates back to the measure $\mu$, finishing the proof.
\endproof

\begin{example}
{When $n = 2$, we may infer immediately from Theorem \ref{thm:cross} that the only Nevanlinna measure whose support is contained in the (closed) first quadrant is the trivial measure.\hfill$\lozenge$}
\end{example}

\section{Properties of measures on the unit polydisk with vanishing mixed Fourier coefficients}\label{sec:measures_polydisk}

In this section, we are going to use the results established for the class of Nevanlinna measures in Section \ref{sec:measures_plane} and translate them to the case of functions from the unit polydisk into the {closed} right half-plane and their associated measures. Even if this is straightforward, we choose to state the properties of the measures explicitely, in order give a complete picture even for the case of the polydisk.

This class of functions is a generalization of Caratheodory functions  and appears at different places, \eg \cite{KoranyiPukanszky1963,McDonald1990,Savchuk2006}. In particular, it is shown in \cite[Theorem 1]{KoranyiPukanszky1963} that they can be characterized by an integral representation in the following sense. A function $f$ maps the polydisk $\D^n$ analytically to the {closed} right half-plane if and only if it admits an integral representation of the form
\begin{equation}\label{eq:intRep_disk_Nvar}
f(\vec{w}) = \I~\Im[f(\vec{0})] + \frac{1}{(2\pi)^n}\int_{[0,2\pi)^n}\left(2\prod_{\ell=1}^n\frac{1}{1-w_\ell\E^{-\I s_\ell}}-1\right)\diff\nu(\vec{s}),
\end{equation}
where $\nu$ is a finite positive Borel measure on $[0,2\pi)^n$ {with vanishing mixed Fourier coefficients}, {\ie} 
\begin{equation}\label{eq:KP_condition}
\int_{[0,2\pi)^n}\E^{\I\:m_1 s_1}\ldots\E^{\I\:m_n s_n}\diff\nu(\vec{s}) = 0
\end{equation}
for any multiindex $\vec{m} \in \Z^n$ with at least one positive entry and at least one negative entry. In \cite{LugerNedic2016,LugerNedic2017}, we have used this characterization and transformed it via a suitable Cayley transform to the poly-upper half-plane. Here, we are going to utilize its inverse transform instead.

To do that in practice, we also need notations for the subsets {of} $[0,2\pi)^n$ that correspond to coordinate-orthogonal subspaces in $\R^n$. A hyperplane of $[0,2\pi)^n$ that is orthogonal to some coordinate axis will be denoted as
$$A_j(p) := \{\vec{s} \in [0,2\pi)^n~|~s_j = p\},$$
where $j \in \{1,2,\ldots,n\}$ and $p \in [0,2\pi)$. Theorem \ref{thm:affine_hyperplane_extraction} and Corollary \ref{coro:affine_hyperplane_extraction} in the case of the polydisk can now be formulated as follows.
 
\begin{thm}\label{thm:measure_affine_subspaces_polydisk}
Let $n \geq 2$, let $\nu$ be a positive Borel measure on $[0,2\pi)^n$ {with vanishing mixed Fourier coefficients} \eqref{eq:KP_condition}. Take an index $j \in \{1,2,\ldots,n\}$ and a point $p \in [0,2\pi)$. Let $\nu|_{A_j(p)}$ denote the restriction of the measure $\nu$ to the hyperplane $A_j(p)$. Then, there {exists} a constant $d_j(p) $ such that
\begin{equation}\label{eq:measure_extracted_form_polydisk}
\nu|_{A_j(p)} = d_j(p)\lambda_{(0,2\pi)^{n{-1}}},
\end{equation}
where {$\lambda_{(0,2\pi)^{n-1}}$} denotes the Lebesgue measure on $(0,2\pi)^{n{-1}}$.
\end{thm}

\begin{remark}
In particular, we note that for a positive Borel measure on $[0,2\pi)^n$ with vanishing mixed Fourier coefficients \eqref{eq:KP_condition}, all points must have zero mass, and, more generally, a statement analogous to Corollary \ref{coro:measure_affine_subspaces} also holds. Furthermore, the area of integration in formula \eqref{eq:KP_condition} may be replaced with the open square $(0,2\pi)^n$.
\end{remark}

\proof{For such a measure $\nu$,} we may use its restriction to the open square {$(0,2\pi)^n$} to build a Nevanlinna measure $\mu$ on $\R^n$ via the mapping $\varphi\colon(0,2\pi) \to \R$, defined by
$$\varphi\colon s \mapsto t:=\I\frac{1+\E^{\I\:s}}{1-\E^{\I\:s}},$$
leading us to define
$$\diff\mu(\vec{t}) := \prod_{j=1}^n{|\varphi'(s_j)|}\diff\nu(\vec{s}).$$
The properties of such a measure $\mu$, as described by Theorem \ref{thm:affine_hyperplane_extraction} and Corollary \ref{coro:affine_hyperplane_extraction}, then translate back to the measure $\nu$ due to the particular way the measure $\mu$ was defined in terms on $\nu$.

This procedure does, in principle, miss a few coordinate parallel affine subspaces of $[0,2\pi)^n$, for example $\{0\} \times [0,2\pi)^{n-1}$, but this is trivially fixed by applying any translation in the definition of the map $\varphi$, that is not an integer multiple of $2\pi$, say
$$\varphi\colon s \mapsto t:=\I\frac{1+\E^{\I\:(s+1)}}{1-\E^{\I\:(s+1)}}.$$
This finishes the proof.
\endproof

\begin{coro}\label{coro:affine_hyperplane_extraction_polydisk}
Let $n \geq 2$, let $f$ be a  function mapping the unit polydisk analytically into the {closed} right half-plane and let $\nu$ be its representing measure {in the sense of representation \eqref{eq:intRep_disk_Nvar}}. Decompose the measure $\nu$ as
$$\nu =  \sum_{i\in I}\nu|_{A_{j_i}(p_i)}+\til{\nu}$$
for some indices $j_i \in \{1,2,\ldots,n\}$ and points {$p_i \in [0,2\pi)$}, where $I {\subseteq \{1,2,\ldots,n\}}$ is a finite set of indices and $\til{\nu}$ the remaining positive Borel measure. Then, the function $f$ can be written as
\begin{equation}\label{eq:intRep_measure_one_extraction_polydisk}
f(\vec{w}) = \sum_{i\in I}d_{j_i}(p_i)\frac{e^{\I\:p_i}+w_{j_i}}{e^{\I\:p_i}-w_{j_i}} + \til{f}(\vec{w}),
\end{equation}
where the function $\til{f}$ {is represented by the measure $\til{\nu}$ in the sense of representation \eqref{eq:intRep_disk_Nvar}}.
\end{coro}

The result on the non-finiteness of Nevanlinna measures, discussed in Proposition \ref{prop:infinite_measure}, for the case of the polydisk may be formulated as follows.

\begin{coro}\label{coro:non_finite_poly}
The function
$$\vec{s} \mapsto \prod_{j=1}^n\frac{1}{s_j^2}$$
is not integrable with respect to any non-trivial positive Borel measure on $[0,2\pi)^n$ {with vanishing mixed Fourier coefficients} \eqref{eq:KP_condition}.
\end{coro}

\proof
Employing the bijection between non-trivial Nevanlinna measures and non-trivial positive Borel measure on $[0,2\pi)^n$ {with vanishing mixed Fourier coefficients} \eqref{eq:KP_condition} as in the proof of Corollary \ref{thm:measure_affine_subspaces_polydisk} yields that
$$\int_{\R^n}\diff\mu(\vec{t}) = \int_{(0,2\pi)^n}\prod_{j=1}^n{|\varphi'(s_j)|}\diff\nu(\vec{s}) = \int_{(0,2\pi)^n}\prod_{j=1}^n\frac{{1}}{1-\cos(s_j)}\diff\nu(\vec{s}) = \infty.$$
Noting that the integrability of the functions $s \mapsto \frac{1}{s^2}$ and $s \mapsto \frac{1}{1-\cos(s)}$ at the point zero is equivalent finishes the proof.
\endproof

In order to  translate the results of Section \ref{subsec:support} to the case of the unit polydisk, we introduce the map $\Phi\colon (0,2\pi)^n \to \R^n$ to be the bijection given as
$$\Phi(\vec{s}) := (\varphi(s_1),\varphi(s_2),\ldots,\varphi(s_n)),$$
where the map $\varphi$ is as in the proof of Corollary \ref{thm:measure_affine_subspaces_polydisk}. 
Under this transformation coordinate-orthogonal hyperplanes {in $\R^n$}  are mapped into coordinate-orthogonal hyperplanes {in $[0,2\pi)^n$}, wheras the image of other affine subspaces are more complicated. The following corollary is, hence, a direct consequence of Theorem \ref{thm:cross}.

\begin{coro}\label{coro:cross_poly}
Let $n \geq 2$, let $\nu$ be a positive Borel measure on $[0,2\pi)^n$ {with vanishing mixed Fourier coefficients} \eqref{eq:KP_condition} and suppose there exist numbers $0 \leq \alpha_j < \beta_j < 2\pi$ for $j = 1,2,\ldots,n$, such that
$$\supp(\nu) \cap \bigcup_{j=1}^n\{\vec{s} \in [0,2\pi)^n~|~\alpha_j < s_j < \beta_j\} = \emptyset.$$
Then, $\nu \equiv 0$.
\end{coro}

\begin{remark}
Due to the fact that the set $[0,2\pi)^n$ is taken as a parametrization of the poly-torus, we could, in Corollary \ref{coro:cross_poly}, just as well consider the union of any combination of sets where $s_j < \alpha_j$ or $\beta_j < s_j$ or both.
\end{remark}

\begin{figure}[!ht]
\includegraphics[width=0.5\linewidth]{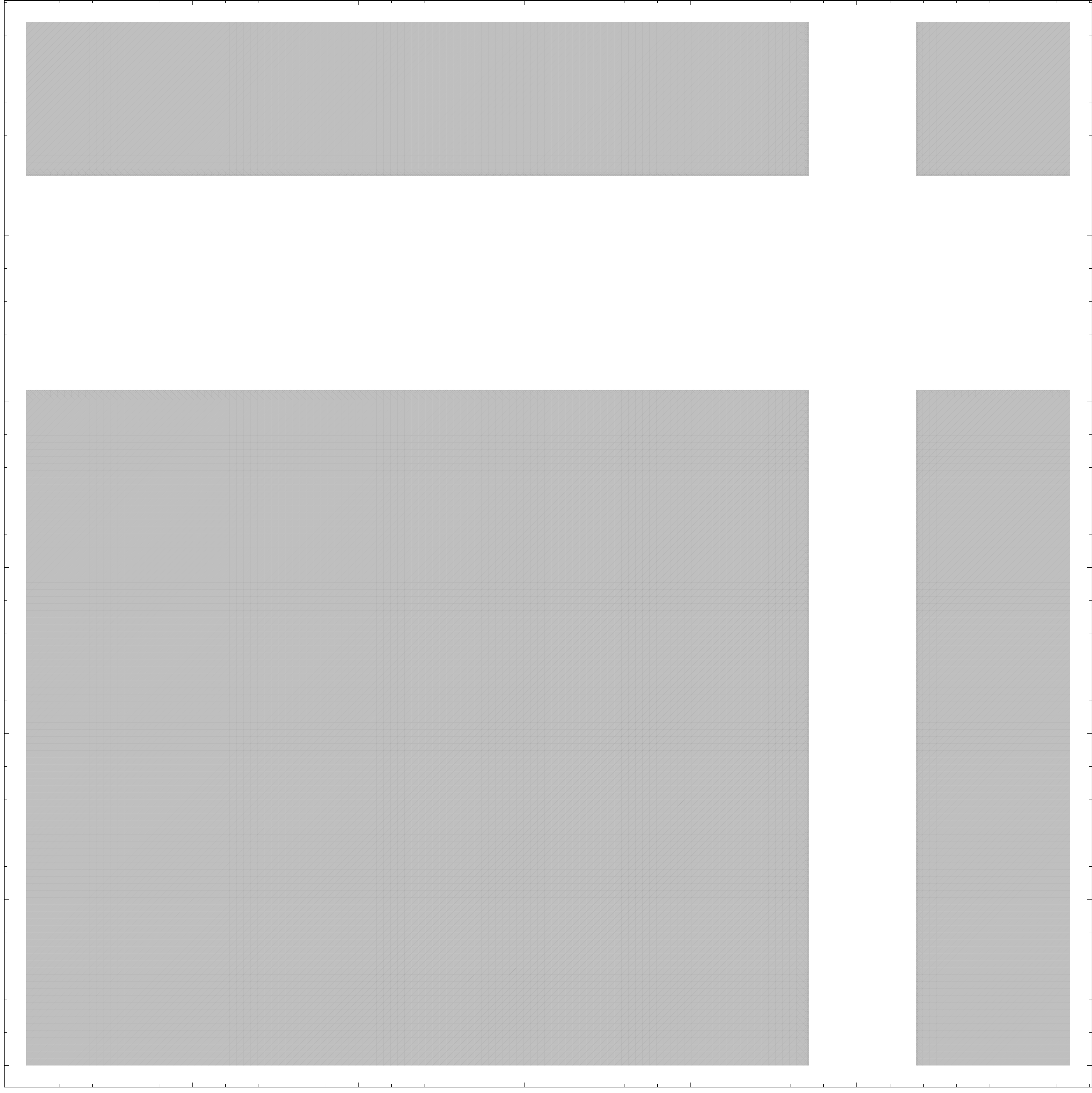}
\caption{The complement of the union of the two strips from Figure \ref{fig:cross}, translated to the case of the polydisk. The plot area is the square $[0,2\pi)^2$.}\label{fig:cross_poly}
\end{figure}

\begin{example}\label{ex:cross_poly}
Let us consider the strips $\{(t_1,t_2) \in \R^2~|~1 < t_1 < 2\}$ and $\{(t_1,t_2) \in \R^2~|~\frac{1}{2} < t_2 < 2\}$ from Figure \ref{fig:cross}. Then, by Theorem \ref{thm:cross}, the only Nevanlinna measure whose support does not intersect the union of these strips in the trivial measure. By Corollary \ref{coro:cross_poly}, we now conclude the only positive Borel measure on $[0,2\pi)^2$ {with vanishing mixed Fourier coefficients} \eqref{eq:KP_condition} whose support does not intersect the union of the strips $\{(s_1,s_2) \in [0,2\pi)^2~|~\varphi^{-1}(1) < s_1 < \varphi^{-1}(2)\}$ and $\{(s_1,s_2) \in [0,2\pi)^2~|~\varphi^{-1}(\tfrac{1}{2}) < s_2 < \varphi^{-1}(2)\}$ is the trivial measure, \cf Figure \ref{fig:cross_poly}. \hfill$\lozenge$
\end{example}

The following corollary is a direct consequence of Theorems \ref{thm:non_par_hyperplanes}, \ref{thm:strips}  and Corollary \ref{coro:minus_inverse}. It will be illustrated with examples below.

\begin{coro}
Let $n \geq 2$ and let $\nu$ be a positive Borel measure on $[0,2\pi)^n$ {with vanishing mixed Fourier coefficients} \eqref{eq:KP_condition} with $\nu|_{A_1(0)} \equiv 0$. Let the sets $\Sigma(\mat{A},\vec{\beta})$ and $S_{j_1,j_2}(\alpha,\beta_1,\beta_2)$ be given as in Theorems \ref{thm:non_par_hyperplanes} and \ref{thm:strips}. Then, it holds that
$$\begin{array}{lcl}
\supp(\nu) \subseteq \Phi^{-1}(\Sigma(\mat{A},\vec{\beta})) & \Longrightarrow & \nu \equiv 0, \\[0.1cm]
\supp(\nu) \subseteq \Phi^{-1}(S_{j_1,j_2}(\alpha,\beta_1,\beta_2)) & \Longrightarrow & \nu \equiv 0, \\[0.1cm]
\supp(\nu) \subseteq \Phi^{-1}(J_1^0(\Sigma(\mat{A},\vec{\beta}) \cap (\R\setminus\{0\})\times\R^{n-1})) & \Longrightarrow & \nu \equiv 0, \\[0.1cm]
\supp(\nu) \subseteq \Phi^{-1}(J_1^0(S_{j_1,j_2}(\alpha,\beta_1,\beta_2) \cap (\R\setminus\{0\})\times\R^{n-1})) & \Longrightarrow & \nu \equiv 0.
\end{array}$$
\end{coro}

\begin{figure}[!ht]
\includegraphics[width=0.5\linewidth]{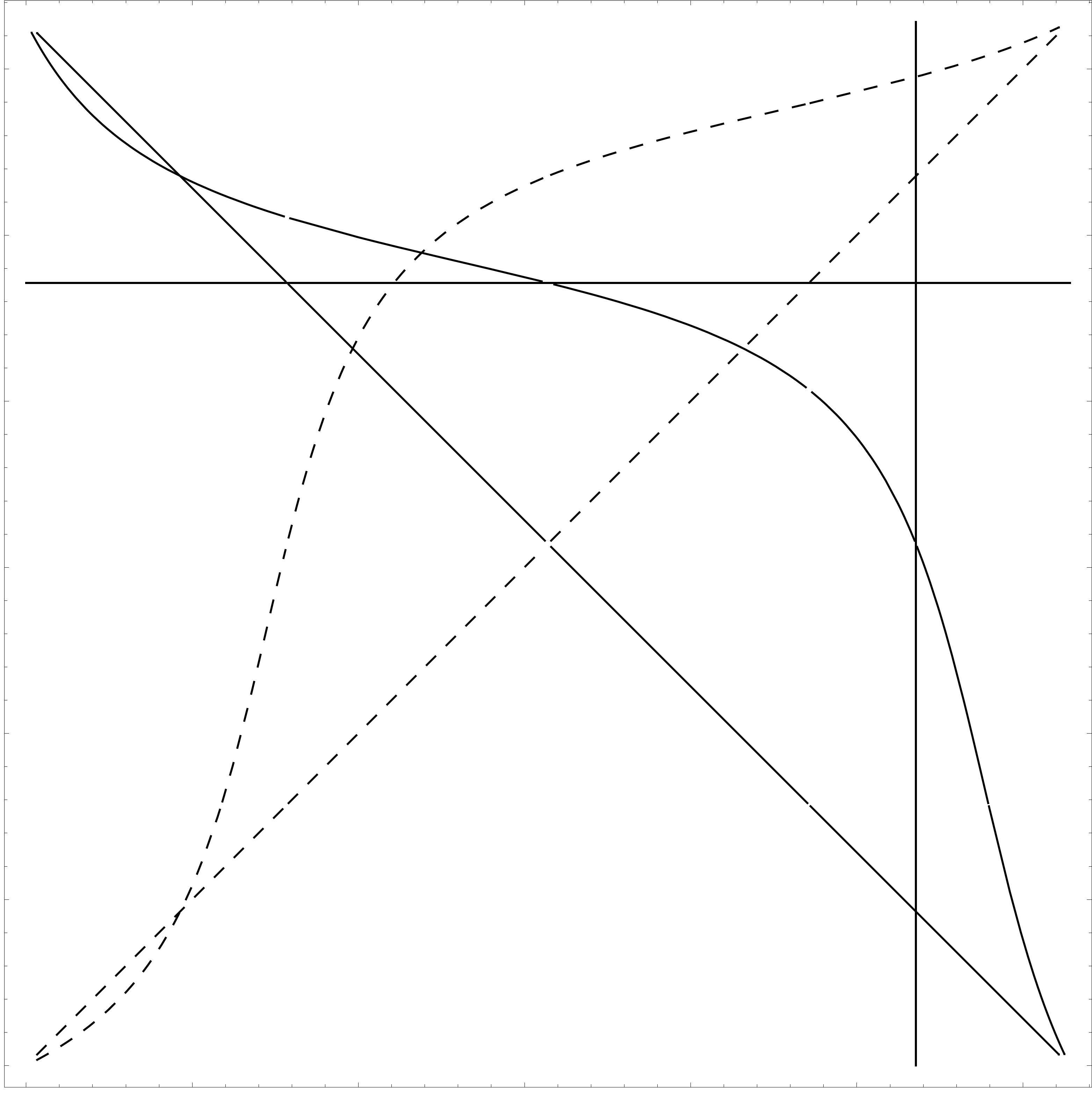}
\caption{The solid lines are those which {can contain} the support of some positive Borel measure on $[0,2\pi)^2$ with vanishing mixed Fourier coefficients \eqref{eq:KP_condition}, while the dashed lines are those {which} cannot contain the support of any such measure, \cf Example \ref{ex:lines_poly}. The plot area is the square $[0,2\pi)^2$.}\label{fig:lines_poly}
\end{figure}

\begin{example}\label{ex:lines_poly}
A non coordinate-orthogonal line in $\R^2$, given by the equation $ t_2 = k\:t_1+ m$ with {$k \in \R\setminus\{0\}$ and $m \in \R$}, yields, on the polydisk side, that
\begin{equation}\label{eq:poly_line}
\varphi(s_2) = k\:\varphi(s_1) + m,
\end{equation}
where $\varphi$ is the same biholomorphism as before. Using the identity
$$\I\frac{1+\E^{\I\:s}}{1-\E^{\I\:s}}=-\cot\left(\frac{s}{2}\right),$$
equation \eqref{eq:poly_line} can be rewritten as $$\cot\left(\frac{s_2}{2}\right) = k\:\cot\left(\frac{s_1}{2}\right) - m,$$ which, in the square $[0,2\pi)^2$, is further equivalent to 
$$s_2 = 2\:{\rm Arccot}\left[k\:\cot\left(\frac{s_1}{2}\right) - m\right].$$
Hence, the curve \eqref{eq:poly_line} is the graph of a function and taking its derivative shows that, for $k>0$, this function is increasing, passing through {the points} $(0,0)$ and $(2\pi, 2\pi)$, whereas, for $k>0$, it is decreasing, passing through {the points} $(0,2\pi)$ and $(2\pi, 0)$. Note that only in the special cases $k=\pm1$ and $m=0$ is the curve actually a straight line. In Figure \ref{fig:lines_poly}, the curves corresponding to the lines from Examples \ref{ex:lines} and \ref{ex:lines_v2} are shown. \hfill$\lozenge$
\end{example}

\begin{figure}[!ht]
\includegraphics[width=0.8\linewidth]{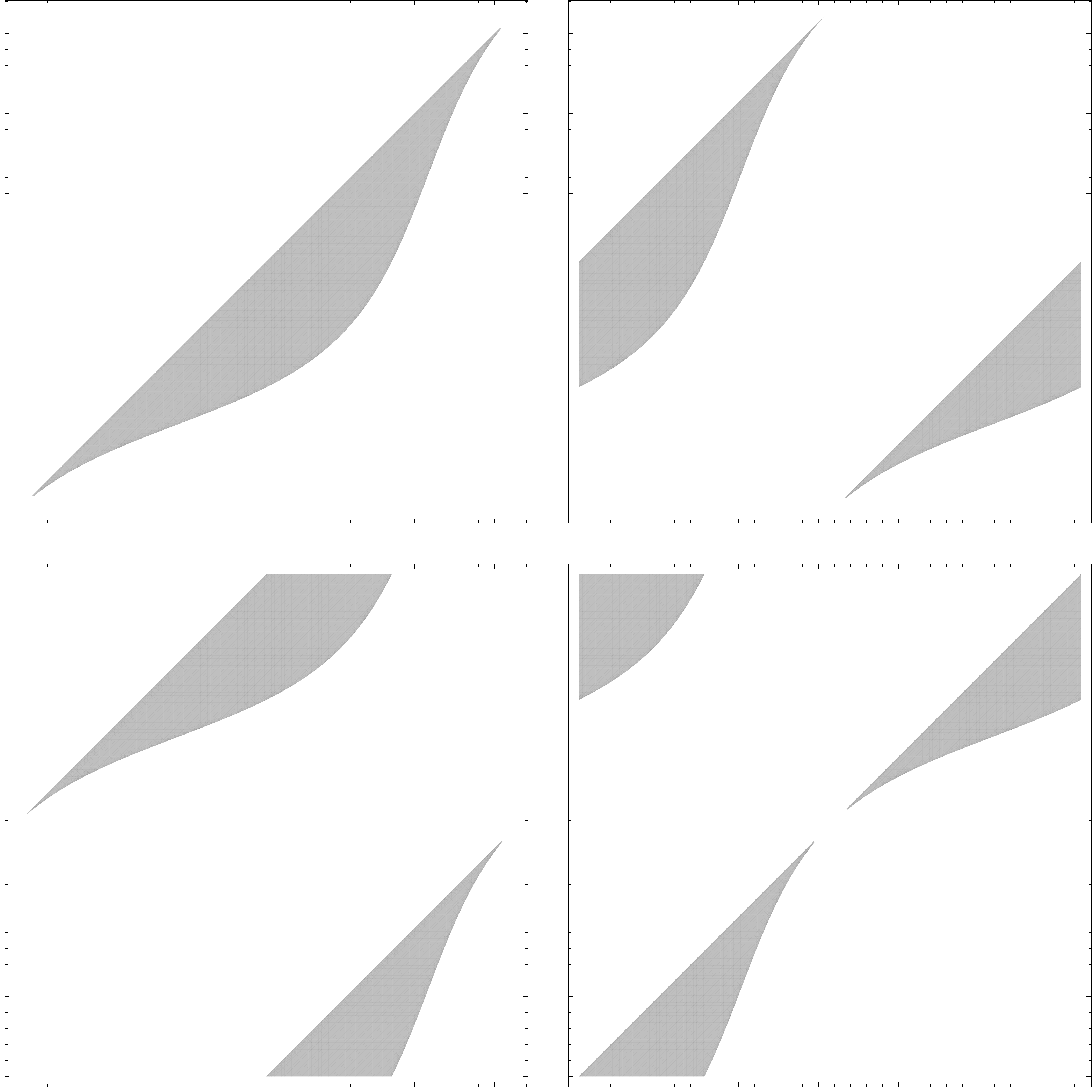}
\caption{The sets $\Phi^{-1}(S_{1,2}(1,-1,0))$ (top left) $\Phi^{-1}(J_1^0(S_{1,2}(1,-1,0) \cap (\R\setminus\{0\}) \times \R))$ (top right), $\Phi^{-1}(J_2^0(S_{1,2}(1,-1,0) \cap \R \times (\R\setminus\{0\})))$ (bottom left) and $\Phi^{-1}(J_1^0J_2^0(S_{1,2}(1,-1,0) \cap (\R\setminus\{0\})^2))$ (bottom right) from Example \ref{ex:strips_poly}. The plot area is always the square $[0,2\pi)^2$.}\label{fig:strips_poly}
\end{figure}

\begin{example}\label{ex:strips_poly}
Let us consider the strip
$$S_{1,2}(1,-1,0) = \{(t_1,t_2) \in \R^2~|~0<t_1-t_2<1\}$$
form Example \ref{ex:strips}, which is bounded by the lines $t_2 = t_1$ and $t_2 = t_1 - 1$ in $\R^2$. Using the information form Example \ref{ex:lines_poly}, we establish that the set $\Phi^{-1}(S_{1,2}(1,-1,0))$ will be bounded by the curves in $[0,2\pi)^2$, given by the equations
$$s_2  = s_1 \quad\text{and}\quad {s_2 = 2\:\mathrm{Arccot} \left[\cot\left(\frac{s_1}{2}\right) + 1\right]}.$$
The boundaries of the set $\Phi^{-1}(J_1^0(S_{1,2}(1,-1,0) \cap (\R\setminus\{0\}) \times \R))$, as well as the sets $\Phi^{-1}(J_2^0(S_{1,2}(1,-1,0) \cap \R \times (\R\setminus\{0\})))$ and $\Phi^{-1}(J_1^0J_2^0(S_{1,2}(1,-1,0) \cap (\R\setminus\{0\})^2))$, may be established analogously, and all four sets are visualized in Figure \ref{fig:strips_poly}. \hfill$\lozenge$
\end{example}

\section*{Acknowledgments}

The authors would like to thank H\r{a}kon Hedenmalm for his inquiries that led us to Theorems \ref{thm:non_par_hyperplanes} and \ref{thm:strips}.

\bibliographystyle{amsplain}
\bibliography{total,total2}

\end{document}